 \documentclass[hidelinks,onefignum,onetabnum]{siamart220329}
 
\usepackage[utf8]{inputenc}
\usepackage{mathtools, amssymb}
\usepackage{enumitem}
\usepackage{xcolor}
\usepackage{kotex}
\usepackage{subfigure}
\usepackage{tabularray}
\usepackage{threeparttable}
\UseTblrLibrary{booktabs}
\usepackage{bbm}
\usepackage{multirow}
\usepackage{stackengine,scalerel}
\usepackage{cleveref}
\usepackage{amsmath}
\usepackage{graphicx}
\usepackage{tikz}
\usetikzlibrary{shapes, arrows.meta, positioning}

\theoremstyle{plain}
\newtheorem{condition}{Condition}

\newcommand{\bi}{\boldsymbol{i}}

\newcommand{\bv}{\boldsymbol{v}}

\newcommand{\bR}{\mathbb{R}}

\newcommand{\cO}{\mathcal{O}}
\newcommand{\cV}{\mathcal{V}}
\newcommand{\cE}{\mathcal{E}}
\newcommand{\cR}{\mathcal{R}}

\newcommand*{\floor}[1]{ \left\lfloor{#1}\right\rfloor }
\newcommand*{\ceil}[1]{ \left\lceil{#1}\right\rceil }
\newcommand{\lrp}[1]{\left({#1}\right)}
\newcommand{\lra}[1]{\left\langle{#1}\right\rangle}

\newcommand{\coe}{\sigma}

\usepackage{lipsum}
\usepackage{amsfonts}
\usepackage{graphicx}
\usepackage{epstopdf}
\usepackage{algorithmic}
\ifpdf
  \DeclareGraphicsExtensions{.eps,.pdf,.png,.jpg}
\else
  \DeclareGraphicsExtensions{.eps}
\fi

\newsiamremark{remark}{Remark}
\newsiamremark{hypothesis}{Hypothesis}
\crefname{hypothesis}{Hypothesis}{Hypotheses}
\newsiamthm{claim}{Claim}

\headers{Localized Estimation of Condition Numbers for MILU}{G. Hwang, Y. Park, Y. Lee, J. Hahn, and M. Kang}

\title{Localized Estimation of Condition Numbers for MILU Preconditioners on a Graph\thanks{Submitted to the editors DATE.
\funding{This work is supported by the National Research Foundation of Korea(NRF) grant funded by the Korea government(MSIT)
 (RS-2024-00421203, RS-2024-00406127, RS-2024-00343226, 2021R1A2C3010887). This project No. 2140/01/01 has received funding from the European Union´s Horizon 2020 research and innovation programme under the Marie Sk\l{}odowska-Curie grant agreement No. 945478.}}}

\author{Geonho Hwang\thanks{Center for AI and Natural Sciences, Korea Institute for Advanced Study
  (\email{hgh2134@kias.re.kr}).}
\and Yesom Park\thanks{Department of Mathematics, University of California 
  (\email{yeisom@math.ucla.edu}).}
  \and Yueun Lee \thanks{Department of Statistics, Seoul National University (\email{mathlife@snu.ac.kr}).}
  \and Jooyoung Hahn \thanks{Department of Mathematics and Descriptive Geometry, Faculty of Civil Engineering, Slovak University of Technology in Bratislava (\email{jooyoung.hahn@stuba.sk}).}
\and Myungjoo Kang\thanks{Department of Mathematical Sciences, Seoul National University (\email{mkang@snu.ac.kr})}}

\usepackage{amsopn}

\ifpdf
\hypersetup{
  pdftitle={Localized Estimation of Condition Numbers for MILU},
  pdfauthor={G. Hwang, Y. Park, Y. Lee, J. Hahn, and M. Kang}
}
\fi

\begin{document}

\maketitle

\begin{abstract}
This paper proposes a theoretical framework for analyzing Modified Incomplete LU (MILU) preconditioners. Considering a generalized MILU preconditioner on a weighted undirected graph with self-loops, we extend its applicability beyond matrices derived by Poisson equation solvers on uniform grids with compact stencils. A major contribution is, a novel measure, the \textit{Localized Estimator of Condition Number (LECN)}, which quantifies the condition number locally at each vertex of the graph. We prove that the maximum value of the LECN provides an upper bound for the condition number of the MILU preconditioned system, offering estimation of the condition number using only local measurements. This localized approach significantly simplifies the condition number estimation and provides a powerful tool or analyzing the MILU preconditioner applied to previously unexplored matrix structures. To demonstrate the usability of LECN analysis, we present three cases: (1) revisit to existing results of MILU preconditioners on uniform grids, (2) analysis of high-order implicit finite difference schemes on wide stencils, and (3) analysis of variable coefficient Poisson equations on hierarchical adaptive grids such as quadtree and octree. For the third case, we also validate LECN analysis numerically on a quadtree.
\end{abstract}

\begin{keywords}
Preconditioner, MILU, Quadtree, High-Order Scheme, Poisson Equation
\end{keywords}

\begin{MSCcodes}
    65F08, 65F35,  05C50, 65N22
\end{MSCcodes}

\section{Introduction}
Iterative methods are widely employed to solve large, sparse linear systems $Ax=b$ of equations, particularly when direct methods are computationally impractical. These methods typically construct a Krylov subspace and project the original system onto it, from which an approximate solution is extracted. The most common iterative techniques include the Conjugate Gradient (CG) method \cite{hestenes1952methods} for symmetric positive-definite systems, and the Generalized Minimal Residual (GMRES) method \cite{saad1986gmres} for non-symmetric systems. Since a larger condition number of $A$ leads to increased computational costs, designing an effective preconditioner to reduce the condition number of the preconditioned system is crucial for enhancing computational performance.

Despite considerable progress in the development of preconditioners \cite{chan1987analysis, kearfott1990preconditioners, axelsson1990algebraic, axelsson1983preconditioning, benzi2004preconditioner}, devising a universally effective preconditioner across diverse problems and matrix structures continues to pose substantial challenges. The effectiveness of a preconditioner largely depends on specific matrix characteristics, such as sparsity patterns and inherent structure. For example, multigrid preconditioners are highly effective to solve matrix systems derived from Poisson-like problems defined on uniform grids. When little else is known about the matrix besides its sparsity pattern, general-purpose or algebraic preconditioners are typically used. These include methods such as Incomplete LU (ILU) \cite{lee2003incomplete, saad1994ilut, malas2007incomplete} factorization, approximate inverse preconditioning \cite{benzi1998sparse, chow1998approximate, chen2001analysis}, and Algebraic Multigrid (AMG) methods \cite{ruge1987algebraic, olson2007algebraic, yang2002boomeramg}. There remains a significant challenge in developing theoretical frameworks to systematically analyze preconditioner performance across various problems. Additionally, many existing techniques need to be examined for their general applicability.

The Modified Incomplete LU (MILU) preconditioner \cite{gustafsson1978class, yoon2017comparison} has shown superior performance compared to other types of ILU preconditioners. One of the main advantages is that the number of nonzero entries in the preconditioned matrix is of the same order of magnitude as that of the original matrix $A$, ensuring that the computational cost of operations remains comparable to that of the original matrix $A$. Moreover, MILU is memory-efficient, requiring only the storage of diagonal entries in addition to the matrix $A$. In \cite{yoon2015analyses}, it is shown that the MILU preconditioner reduced the condition number from an order of $\cO\left(h^{-2}\right)$ to $\cO\left(h^{-1}\right)$, where $h$ is the grid size, for the finite difference scheme in \cite{gibou2002second} to solve the two dimensional (2D) Poisson equation with Dirichlet boundary conditions. In \cite{hwang2024analysis}, this result is extended to three dimension (3D) and a variant, Sectored-MILU (SMILU), is designed to also reduce the condition number to a comparable order. Experimental results numerically verified that SMILU provides more reduction of the condition number than MILU. The analysis of the mentioned preconditioners was limited to uniform grids and Poisson equations.

In this paper, we propose a theoretical framework for analyzing MILU preconditioners in a broader range of problems. When the coefficient matrix $A$ is represented by a weighted adjacency matrix of a weighted undirected graph with self-loops, a general MILU on the graph extends the applicability of MILU beyond previously studied numerical methods on uniform grids with compact stencils. The core idea of the proposed theory is, a novel measure, \textit{Localized Estimator of Condition Number (LECN)}, which quantifies the condition number locally at each vertex of the graph. We prove that the maximum value of the LECN serves as an upper bound for the condition number of the MILU preconditioned system. This framework, \textit{LECN analysis}, provides a method to approximate the condition number by local characteristics at a vertex of the associated graph, offering an efficient and scalable method for analyzing the condition number. 

To demonstrate the usability of LECN analysis, we present three cases. First, LECN analysis can provide more refined proofs than~\cite{yoon2015analyses, hwang2024analysis}. Second, we extend the analysis to investigate high-order implicit finite difference schemes on wide stencils~\cite{claerbout1985craft, zapata2017high}, proving that MILU achieves a condition number of order $\cO\left(h^{-1}\right)$. Third, we extend the analysis to hierarchical adaptive grids, such as quadtree and octree grids, to estimate the condition number of MILU preconditioned system from finite volume method \cite{frolkovivc2021flux} to solve variable coefficient Poisson equations. Since sizes of cells and neighborhood structures are changed locally on quadtree and octree grids, it mathematically challenges to analyze the condition number. We prove that MILU reduces the condition number to the order of $\cO\left(\bar{h}^{-1}\right)$, where $\bar{h}$ is the size of the smallest cell in the grid. This case is particularly significant, as LECN analysis is not only applicable to analyze variable coefficient Poisson equations but also capable to analyze linear systems on hierarchical adaptive grids. We also numerically verify theoretical results for the last challenging case. These cases demonstrate how a localized approach of LECN analysis faciliates to study MILU preconditioners on weighted undirected graphs with self-loops. This allows for the exploration of high-order schemes and hierarchical adaptive grids, which previous literature has not theoretically covered.

The rest of the paper is organized as follows. In \cref{sec:MILU-type Preconditioner}, the mathematical notation and the scope of problems are introduced and then we present the LECN analysis. In \cref{sec:Poissson_uniform_grid}, the established results of MILU preconditioners are revisited by the LECN analysis and we show that LECN analysis facilitates the analysis of high-order schemes. In \cref{sec:MILU_variable_quad_octree}, we provide an analysis of condition number for matrix systems arising from variable coefficient Poisson equations on hierarchical adaptive grids such as quadtree and octree.

\section{Localized Estimator of Condition Number}\label{sec:MILU-type Preconditioner}

\subsection{Problem Setting}\label{sec:notation_defintion}
Consider a linear system of equations 
\begin{equation}\label{eq:linear_sys}
	Ax = y,
\end{equation}
where the coefficient $A \in \mathbb{R}^{N \times N}$ is a symmetric positive definite (SPD) M-matrix:
\begin{equation}\label{eq:def_A}
	A_{K,K'} = \begin{cases}
		- c_{K,K'}&\text{ if } K\neq K',
		\\ \sum_{K'\neq K} c_{K,K'} + b_K &\text{ if } K=K',
	\end{cases} 
\end{equation}
where $c_{K,K'} = c_{K',K} \geq 0$ and $b_K\geq 0$ for $K$, $K' \in \{ 1, \ldots, N\}$. To ensure the positive definiteness of $A$, at least one of $b_K$ must be strictly positive.
The condition number of the coefficient matrix $A$ is defined
\begin{equation*}
	\kappa\left(A\right)=\frac{\lambda_{\max}}{\lambda_{\min}},
\end{equation*}
where $\lambda_{\max}$ and $\lambda_{\min}$ are the maximal and minimal eigenvalues of $A$, respectively. It is a crucial measure that determines the convergence rate of iterative methods such as the conjugate gradient method (CG) \cite{hestenes1952methods} or the generalized minimal residual method (GMRES) \cite{saad1986gmres}. A large condition number slows down the convergence, requiring more computational cost. A preconditioner, a matrix $M\in \mathbb{R}^{N \times N}$, is designed to improve the condition number and thus accelerate the convergence of iterative solvers. A preconditioned system is given by
\begin{equation*}
	M^{-1}Ax = M^{-1}y,
\end{equation*}
where the preconditioned matrix $M^{-1}A$ ideally has a smaller condition number than $A$, closer to $1$.

To achieve a theoretical framework of analyzing MILU preconditioners, we reinterpret the coefficient matrix $A \in \mathbb{R}^{N \times N}$ \eqref{eq:def_A} by a weighted adjacency matrix of a weighted undirected graph with self-loops $G = \lrp{\mathcal{V}, \mathcal{E}, w}$. The set of vertices in the graph is $\mathcal{V} = \{ 1, \dots, N\}$ and the set of edges is $\mathcal{E} = \{e_{K, K'} : A_{K, K'} \neq 0, \: K, \: K' \in \mathcal{V} \}$. Note that $e_{K, K'} = e_{K', K}$. The weight function $w:E \rightarrow \mathbb{R}$ assigns a weight to an edge: For $K$ and $K' \in \mathcal{V}$ and $e_{K,K'} \in \mathcal{E}$,
\begin{align*}
	w(e_{K,K'}) = A_{K,K'}.
\end{align*}
The sparsity of the matrix $A$ directly corresponds to the connectivity of the edges in the graph. For a vertex $K \in \mathcal{V}$, we denote $n_0(K)$ as the set of neighboring vertices of $K$ connected by an edge but not the self-loop at $K$:
\begin{equation*}
	n_0(K) = \{K'\in \mathcal{V}:e_{K,K'} \in \mathcal{E}, K \neq K'\}.
\end{equation*}
Letting a vector space $S = \mathbb{R}^{\mathcal{V}}$, the set of matrices that map from $S$ to $S$ is denoted by $\mathcal{M}(S)$. In \cref{sec:Poissson_uniform_grid} and \cref{sec:MILU_variable_quad_octree}, we provide more specific examples of the coefficient matrix $A$ above.

\subsection{MILU on Graphs}\label{sec:milu_graph}
To define a MILU preconditioner on a weighted undirected graph with self-loops $G$ from the coefficient matrix $A$~\eqref{eq:def_A}, we first introduce a partial order on the graph $(G, \prec)$ and then the MILU preconditioner is recursively calculated by a specific ordering. 
\begin{definition}[Order on Graph] For a graph $G$, a \textbf{partial order} on $G$ is a strict partial order $\prec$ defined on the vertices of $G$, where a strict order is assigned between any two neighboring vertices. Specifically, for a vertex $K_1\in\cV$ and its neighboring vertex $K_2\in n_0\left(K_1\right)$, either
\begin{equation*}
	K_1\prec K_2 \text{ or } K_2\prec K_1.
\end{equation*}
Since any partial order can be extended to a total order, the partial order $\prec$ assigns a total order to the entire graph $G$.
\end{definition}
\begin{definition}[Precursor and Successor] For each vertex $K\in\cV$, the \textbf{precursor} of $K$ is defined by the set of neighboring vertices with a lower order than $K$:
\begin{equation*}
	p(K)= \left\{K'\in n_0(K): K'\prec K\right\},
\end{equation*}
while the set of neighboring vertices with a higher order is termed the \textbf{successor}:
\begin{equation*}
	s(K)= \left\{K'\in n_0(K): K\prec K'\right\}.
\end{equation*}
\end{definition}

Note that $n_0(K)=p(K)\sqcup s(K)$. With $(G, \prec)$, the matrix $A$ \eqref{eq:def_A} can be partitioned into $A = L+D + L^\top$, where the (strict) lower triangular part $L$ contains the connections corresponding to the precursors, the (strict) upper triangular part $L^\top$ contains the connections corresponding to the successors, and the diagonal $D$ represents the weights on the self-loops.

\begin{definition}[MILU on Graph]
    Given a weighted undirected graph with self-loops $G=\left(\cV,\cE, w\right)$ with a partial order $\prec$, the \textbf{MILU preconditioner} for the matrix $A=L+D+L^\top \in \mathcal{M}\left(S\right)$ is defined as follows: 
    \begin{equation}\label{eq:milu}
    M = (L+ E)E^{-1}(L+E)^{\top},
\end{equation}
where $E\in \mathcal{M}\left(S\right)$ is a diagonal matrix, defined such that the row sums of 
\begin{equation*}
    R \coloneqq M - A=D  - E  - LE^{-1}L^\top
\end{equation*}
are zero. 
Specifically, the diagonal entries $e_K = E_{K,K}$ of $E$ is defined recursively according to the order of the vertices as follows:
\begin{align}\label{eq:e_recursive_equation}
\begin{split}
    e_{K} & = A_{K,K} - \sum_{K_1\in p(K), K_2\in s(K_1)}\frac{A_{K,K_1}A_{K_1,K_2}}{e_{K_1}}\\
    & = b_K - \sum_{K_1\in n_0(K)} A_{K,K_1} - \sum_{K_1\in p(K), K_2\in s(K_1)}\frac{A_{K,K_1}A_{K_1,K_2}}{e_{K_1}}.
     \end{split}
\end{align}
\end{definition}
Note that if $b_K\ge0$ for all $K\in \mathcal{V}$, then $e_k$ is always well-defined. This is evident through the following argument: For $K\in \mathcal{V}$ such that $K_1\nprec K$ for any $K_1\in \mathcal{V}$, $e_K = A_{K,K}\ge b_K - \sum_{K_1\in s(K)}A_{K,K_1}$.
By mathematical induction, assume that $e_{K_1}\ge \sum_{K_2\in s(K_1)}-A_{K_1,K_2}$ for any $K_1\prec K$, then 
\begin{equation*}
    e_K\ge  b_K - \sum_{K_1\in n_0(K)} A_{K,K_1} + \sum_{K_1\in p(K)} A_{K,K_1} \ge \sum_{K_1\in s(K)}-A_{K,K_1}.
\end{equation*}
Thus, the induction hypothesis holds for all $K\in\mathcal{V}$, hence $e_K\ge \sum_{K_1\in s(K)}-A_{K,K_1}$ for all $K\in \mathcal{V}$, which ensures the well-definedness of $e_K$.

We emphasize the flexibility in defining different orderings for a fixed graph, where the weighted adjacency matrix $A$, as defined in \eqref{eq:def_A}, remains unchanged. This variability in ordering enables the derivation of distinct MILU preconditioners from a single SPD M-matrix 
$A$, depending on the chosen graph ordering.
Consequently, the condition number of the system preconditioned by the MILU is influenced by the graph's ordering. The standard MILU preconditioner \cite{gustafsson1978class,yoon2015analyses} and the Sectored MILU (SMILU) \cite{hwang2024analysis} are derived by two different orders, lexicographic and quadruplicate/octuplicate orders, respectively, on the uniform grid. In \cref{sec:validate_poisson}, we show a theoretical insight to explain which order on the same graph can lead to more efficient computation based on the LECN analysis tool in the following section.


\subsection{LECN Analysis for MILU}\label{sec:lecn_analysis}
We develop a theoretical framework to analyze the condition number of the MILU preconditioned system. Our method employs a local measure at each vertex, quantifying the ratio of the diagonal entries of $E$ to the entries associated with its successors. We define this ratio as the \textit{Localized Estimator of Condition Number (LECN)}.

\begin{definition}[Localized Estimator of Condition Number (LECN)]
Considering a partial order on a weighted undirected graph with self-loops $G=\left(\cV,\cE, w\right)$, its weighted adjacent matrix $A\in \mathcal{M}\left(S\right)$, and MILU preconditioner with $e_K$ in \eqref{eq:e_recursive_equation}, \textbf{Localized Estimator of Condition Number (LECN)} $\tau_{K}$ is defined at $K\in \mathcal{V}$: 
\begin{equation}\label{eq:tau}
    \tau_{K} = \frac{e_K}{e_K + \sum_{K_1\in s(K)}A_{K,K_1}}.
\end{equation}
\end{definition}

The primary result of this paper, stated in the following proposition, proves that LECN provides an upper bound for the condition number of the MILU preconditioned system.
\begin{proposition}[LECN analysis]\label{proposition:tau_e}
     For a matrix A \eqref{eq:def_A}, MILU preconditioner $M$ \eqref{eq:milu}, and $\tau_K$ \eqref{eq:tau} for each $K\in\cV$, the following inequality holds:
\begin{equation*}
    \kappa\lrp{M^{-1}A} \le \max_{K\in \mathcal{V}}\tau_K.
\end{equation*}
\end{proposition}
\begin{proof}
As $R$ is row-sum zero, for $\bv\in S$, the following equation holds:
\begin{equation}\label{eq:pf_prop_tau_e_eq_1}
    \lra{-R\bv,\bv} = \sum_{K\in \mathcal{V}}\sum_{\substack{K_1, K_2\in s(K),\\ K_1<K_2}} \frac{A_{K,K_1}A_{K,K_2}}{e_K}\lrp{v_{K_1}- v_{K_2}}^2 .
\end{equation}
We further attain that
  \begin{align}\label{eq:pf_prop_tau_e_eq_2}
  \begin{split}
          &\sum_{K\in \mathcal{V}}\sum_{\substack{K_1, K_2\in s(K),\\ K_1<K_2}}\frac{A_{K,K_1}A_{K,K_2}}{e_K}\lrp{v_{K_1}- v_{K_2}}^2 \\
         & \le \sum_{K\in \mathcal{V}}\lrp{ \lrp{\frac{\sum_{K_2\in s(K)}A_{K,K_2}}{e_K}}\sum_{K_1\in s(K)}A_{K,K_1}\lrp{v_{K}- v_{K_1}}^2},
          \end{split}
  \end{align}
  where the inequality holds by the following lemma:
  \begin{lemma}
      For $n,i\in \mathbb{N}$ satisfying $i\le n$ and $a_i\in \mathbb{R}$, the following inequality holds for arbitrary $x_i\in \mathbb{R}$ and $x\in \mathbb{R}$:
      \begin{equation*}
          \sum_{1\le i<j\le n}a_ia_j(x_i-x_j)^2\le \lrp{\sum_{i=1}^na_i}\sum_{j=1}^na_j(x_j-x)^2.
      \end{equation*}
  \end{lemma}
  \begin{proof}
      The inequality is equivalent to 
      \begin{equation*}
          0\le \lrp{\sum_{j=1}^na_j(x_j-x)}^2.
      \end{equation*}
     which completes the proof.
  \end{proof}

  From \eqref{eq:pf_prop_tau_e_eq_1} and \eqref{eq:pf_prop_tau_e_eq_2}, we have
  \begin{align*}
       \lra{-R\bv, \bv} &= \sum_{K\in \mathcal{V}}\sum_{\substack{K_1, K_2\in s(K),\\ K_1<K_2}}\frac{A_{K,K_1}A_{K,K_2}}{e_K}\lrp{v_{K_1}- v_{K_2}}^2 
        \\ & \le \max_{K\in \mathcal{V}} \lrp{\frac{\sum_{K_2\in s(K)}-A_{K,K_2}}{e_K}} \sum_{K\in \mathcal{V}} \sum_{K_1\in s(K)}-A_{K,K_1}\lrp{v_{K}- v_{K_1}}^2\\
        & \le  \max_{K\in \mathcal{V}} \lrp{\frac{\sum_{K_2\in s(K)}-A_{K,K_2}}{e_K}}  \lra{A\bv, \bv}.
  \end{align*}
Therefore, we have
 \begin{align*}
   1\le   \frac{\lra{\bv, A\bv}}{\lra{\bv, M\bv}} = \frac{1}{1 - \frac{\lra{\bv, -R\bv}}{\lra{\bv, A\bv}}} & \le \frac{1}{1 - \max_{K\in \mathcal{V}} \lrp{\frac{\sum_{K_1\in s(K)}-A_{K,K_1}}{e_K}}} 
    \\ &= \max_{K\in \mathcal{V}} \frac{e_K}{e_K + \sum_{K_1\in s(K)}A_{K,K_1}}
   = \max_{K\in \mathcal{V}} \tau_K ,
\end{align*}
where the first inequality holds as $ \lra{-R\bv, \bv}\ge 0$, as $e_K\ge 0$.
Together with Rayleigh quotient, this completes the proof.
\end{proof}
%

The main strength of the LECN analysis is the estimation of the condition number through localized calculations that require only the consideration of each vertex's neighborhood connected by an edge. This localized approach is useful for handling more complex configurations, such as wide stencils or hierarchical adaptive grids, which we demonstrate in \cref{sec:high_order} and \cref{sec:MILU_variable_quad_octree}. According to the LECN analysis, it is essential to monitor the value of $\tau_K$ and we provide a useful lemma for calculating $\tau_K$.

  \begin{lemma}\label{lemma:tau_evolove}
  For $K\in\cV$,
      \begin{equation*}
    \tau_K = 1 +\frac{\sum_{K_2\in s(K)}-A_{K,K_2}}{b_K + \sum_{K_1\in p(K)} \frac{-A_{K,K_1}}{\tau_{K_1}}}.
\end{equation*}
  \end{lemma}
\begin{proof}
    Recall that 
\begin{equation*}
   e_K  = b_K - \sum_{K_1\in n_0(K)} A_{K,K_1} - \sum_{K_1\in p(K), K_2\in s(K_1)}\frac{A_{K,K_1}A_{K_1,K_2}}{e_{K_1}}.
\end{equation*}
It follows that
 \begin{align*}
  \begin{split}
    e_K + \sum_{K_1\in s(K)}A_{K,K_1} & = b_K - \sum_{K_1\in p(K)} A_{K,K_1} - \sum_{K_1\in p(K), K_2\in s(K_1)}\frac{A_{K,K_1}A_{K_1,K_2}}{e_{K_1}}
    \\ 
    &= b_K - \sum_{K_1\in p(K)} A_{K,K_1}\lrp{1 + \frac{\sum_{K_2\in s(K_1)}A_{K_1,K_2}}{e_{K_1}}}.
\end{split}
\end{align*}
Therefore, we can attain the following identity
\begin{equation*}
    \frac{\sum_{K_1\in s(K)}-A_{K,K_1}}{\tau_K - 1} =  b_K - \sum_{K_1\in p(K)} \frac{A_{K,K_1}}{\tau_{K_1}},
\end{equation*}
which leads to
  \begin{equation*}
    \tau_K = 1 +\frac{\sum_{K_1\in s(K)}-A_{K,K_1}}{b_K + \sum_{K_1\in p(K)} \frac{-A_{K,K_1}}{\tau_{K_1}}}.
\end{equation*}
\end{proof}

\begin{remark}\label{rmk:Obs_LECN_anal}
A key insight from the LECN analysis, supported by \cref{lemma:tau_evolove}, is that at each vertex $K$, it is less desirable to have a larger sum of matrix values for successors, $\sum_{K_2\in s(K)}-A_{K,K_2}$, compared to precursors, $\sum_{K_1\in p(K)}-A_{K,K_1}$, which causes an increase of the value $\tau_K$. In particular, if vertices $K$ with $b_K=0$ have no precursors, $\tau_K$ becomes infinite. Therefore, such vertices must have at least one precursor. In this regard, \cref{lemma:tau_evolove} provides guidance on how to appropriately order the vertices of $A$, taking into account both the values of successors and precursors.
\end{remark}

\section{Poisson Equation on Uniform Grids}\label{sec:Poissson_uniform_grid}

In this section, we apply LECN analysis to estimate the condition number of MILU preconditioned systems arising from finite difference methods (FDM) based on uniform grids for solving the Poisson equation with Dirichlet boundary condition on a bounded domain $\Omega\in\bR^d$ for $d=2,3$:
\begin{equation}\label{eq:poisson}
	\begin{cases}
		- \triangle u(x) = f(x), &x\in \Omega,
		\\ u(x) = g(x), &x\in \partial\Omega,
	\end{cases}
\end{equation}
where $f:\Omega\rightarrow \bR$ is a source function and $g:\partial\Omega\rightarrow\bR$ is a boundary function.
In \cref{sec:validate_poisson}, we revisit existing MILU analyses \cite{yoon2015analyses, hwang2024analysis} for the second-order FDM \cite{gibou2002second} in the view of LECN analysis, leading to more refined proofs. In \cref{sec:high_order}, we use the LECN analysis to analyze high-order FDMs \cite{claerbout1985craft, zapata2017high} on wide stencils, highlighting the effectiveness of LECN in handling complex graph structures.

\subsection{Revisit of Existing Analysis}\label{sec:validate_poisson}
We employ the LECN analysis to revisit the theoretical analyses of MILU and SMILU for the second-order FDM \cite{gibou2002second} for solving \eqref{eq:poisson}.
The computational domain is discretized using a uniform grid of grid size $h$, with the grid points within the irregular domain $\Omega$ forming the vertices of the graph. For each vertex $K\in\cV$, its neighbors $n_0\lrp{K}$ are defined as the grid points $K'\in\cV$ directly connected to $K$ along grid lines. For a vertex $K\in\cV$, the set of grid points directly connected to $K$ along a grid line that lie outside the domain are denoted as $B\lrp{K}$. Then, the coefficient matrix derived by the discretization in~\cite{gibou2002second} is given by the following SPD M-matrix:
\begin{equation}\label{eq:gibou}
\small
A_{K_1, K_2} = \begin{cases}
		- h^{-2} & \text{ if $K_1\in n_0\lrp{K_2}$ or $K_2\in n_0\lrp{K_1}$},\\  
		|n_0(K)|h^{-2} + \sum_{K' \in B(K)}\left(h_{K,K'}h\right)^{-1} & \text{ if $K_1 = K_2 (=K)$},\\
		0  & \text{ otherwise},
	\end{cases}
\end{equation}
where $h_{K,K'}$ is a distance from $K$ to $\partial \Omega$ along the grid line to the grid point $K'\in B\lrp{K}$. Note that the condition number of $A$ is the order $\cO\left(h^{-2}\right)$ \cite{yoon2015analyses}.

The standard MILU preconditioner, defined by a lexicographical order on the weighted undirected graph with self-loops from the matrix~\eqref{eq:gibou}, has been proven to reduce the condition number to order $\cO\lrp{h^{-1}}$ in 2D \cite{yoon2015analyses} and 3D \cite{hwang2024analysis} domains. The Sectored-MILU (SMILU) \cite{hwang2024analysis}, defined in 2D and 3D by quadruplicate and octuplicate orders on the same graph, respectively, also achieves a similar reduction of the condition number. In the following \cref{cor:standard_milu} and~\cref{cor:smilu}, we provide more refined proofs directly from \cref{proposition:tau_e} than~\cite{yoon2015analyses, hwang2024analysis}. Additionally, we explain theoretical insights into the superior performance of SMILU over MILU.

\begin{corollary}\label{cor:standard_milu}
	For the coefficient matrix $A$ defined in \eqref{eq:gibou} with grid size $h$ and the MILU preconditioner $M$ defined on the lexicographical order, the condition number of preconditioned system is bounded by
	\begin{equation*}
		\kappa\lrp{M^{-1}A} \le 1+d+\frac{d\ell_{\max}}{h},
	\end{equation*}
	where $\ell_{\max}$ is the length of the longest side among the axes of the computational domain.
\end{corollary}
\begin{proof} 
	In the lexicographical order, for an interior grid $K$ where $B\left(K\right)=\emptyset$, both precursors and successors are assigned
	$d$ elements each.
	Therefore, for an interior grid $K$,
	\begin{align}\label{eq:pf_cor_milu_poisson_tau_flow}
		\begin{split}
			\tau_K &= 1 +\frac{\sum_{K_1\in s(K)}-A_{K,K_1}}{\sum_{K_1\in p(K)} \frac{-A_{K,K_1}}{\tau_{K_1}}}
			= 1 +\frac{\sum_{K_1\in s(K)}\frac{1}{h^2}}{\sum_{K_1\in p(K)} \frac{1}{h^2}\frac{1}{\tau_{K_1}}} \\
			& = 1 +\frac{|s(K)|}{\sum_{K_1\in p(K)}\frac{1}{\tau_{K_1}}} 
			= 1 +\frac{|p(K)|}{\sum_{K_1\in p(K)}\frac{1}{\tau_{K_1}}}
			\le 1 + \max_{K_1\in p(K)}\tau_{K_1}
			.
		\end{split}
	\end{align}
	For each boundary grid $K$, that is, $B\left(K\right) \neq \emptyset$
	\begin{align}\label{eq:pf_cor_milu_poisson_tau_boundary}
		\begin{split}
			\tau_K & = 1 +\frac{\sum_{K_1\in s(K)}-A_{K,K_1}}{\sum_{K'\in B\lrp{K}}  \frac{1}{hh_{K,K'}} + \sum_{K_1\in p(K)} \frac{-A_{K,K_1}}{\tau_{K_1}}}
			\le 1 +\frac{\sum_{K_1\in s(K)}\frac{1}{h^2}}{\sum_{K'\in B\lrp{K}}  \frac{1}{hh_{K,K'}}}\\
			& \leq 1 +\frac{\sum_{K_1\in s(K)}\frac{1}{h^2}}{\sum_{K'\in B\lrp{K}}  \frac{1}{h^2}}
			\leq 1 + |s(K)| 
			\le 1 + d,
		\end{split}
	\end{align}
	where the second inequality is induced from the fact that $h_{K,K'}<h$ for $K'\in B\lrp{K}$.
	
	Given that the lexicographical order is defined over the entire graph, each interior grid $K$ is connected to a boundary node through its precursors inductively. Applying \eqref{eq:pf_cor_milu_poisson_tau_flow} and \eqref{eq:pf_cor_milu_poisson_tau_boundary}, we obtain the following estimate for every non-boundary node $K$ 
	\begin{equation*}
		\tau_K\le 1+d+l_K,
	\end{equation*}
	where $l_K$ represents the increasing path from a boundary note to node $K$. If $K$ is the boundary node, then $l_k=0$. 
	Therefore, we arrive at
	\begin{equation*}
		\tau_K\le 1+d+\max_K l_K,
	\end{equation*}
	where the maximum length $\max_K l_K\leq\frac{d\ell_{\max}}{h}$ in the case of the lexicographical order.
	In conjunction with \cref{proposition:tau_e}, it concludes the proof.
\end{proof}
\begin{corollary}\label{cor:smilu}
	For coefficient matrix $A$ defined in \eqref{eq:gibou} with grid size $h$ and the 
	SMILU preconditioner $M$ defined on the quadruplicate/octuplicate order, the condition number of preconditioned system is bounded by
	\begin{equation*}
		\kappa\lrp{M^{-1}A} \le 1+d+\frac{d\ell_{\max}}{2h},
	\end{equation*}
	where $\ell_{\max}$ is the length of the longest side among the axes of the computational domain.
\end{corollary}
\begin{proof}
	SMILU defines the MILU preconditioner by dividing the computational domain into sectors along each axis of the dimensions, applying a lexicographical ordering starting from the corners of each sector. This ordering ensures that all nodes are drawn towards the center of the domain, guaranteeing that the number of successors is always less than the number of precursors, that is, $\left\vert s\left(K\right)\right\vert\leq \left\vert p\left(K\right)\right\vert$ for every node $K$. By applying the proof method of \cref{cor:standard_milu} to each sector, the similar result is obtained as in \cref{cor:standard_milu} as follows:
	\begin{equation*}
		\tau_K\le 1+d+\max_K l_K,
	\end{equation*}
	Since the domain is divided into halves along each axis, the maximum length of the path from a node $K$ to the boundary is reduced by a factor of two, resulting in $\max_K l_K\leq\frac{d\ell_{\max}}{2}$.
\end{proof}

\Cref{cor:standard_milu} and \cref{cor:smilu} show that LECN analysis provides a theoretical insight of the numerical observation that the quadruplicate/octuplicate order used to define SMILU results in a smaller condition number than the lexicographical order used to define the standard MILU. 
For a vertex $K\in\cV$ where $B\lrp{K}=\emptyset$, the weight of all neighboring vertices $K_1\in n_0\lrp{K}$ is uniformly $A_{K,K_1}=-h^{-2}$, which implies that $\sum_{K_1\in p(K)}-A_{K,K_1}= \left\vert p\lrp{K}\right\vert h^{-2}$ and  $\sum_{K_1\in s(K)}-A_{K,K_1}= \left\vert s\lrp{K}\right\vert h^{-2}$.
Therefore, as noted in~\cref{rmk:Obs_LECN_anal}, having fewer successors than precursors is beneficial for reducing the LECN value $\tau_K$.
In the lexicographical order, $\left\vert p\lrp{K}\right\vert=\left\vert s\lrp{K}\right\vert$, whereas in a quadruplicate/octuplicate order, $\left\vert p\lrp{K}\right\vert\geq\left\vert s\lrp{K}\right\vert$.
This difference results in SMILU having a smaller LECN value than MILU, as demonstrated in \cref{cor:standard_milu} and \cref{cor:smilu}, which would ultimately contribute to the improved performance of SMILU.

\subsection{LECN analysis for HIFD}\label{sec:high_order}

\begin{figure}[htbp]
\begin{center}
\begin{tabular}{ccc}
\includegraphics[height=0.15\textwidth]{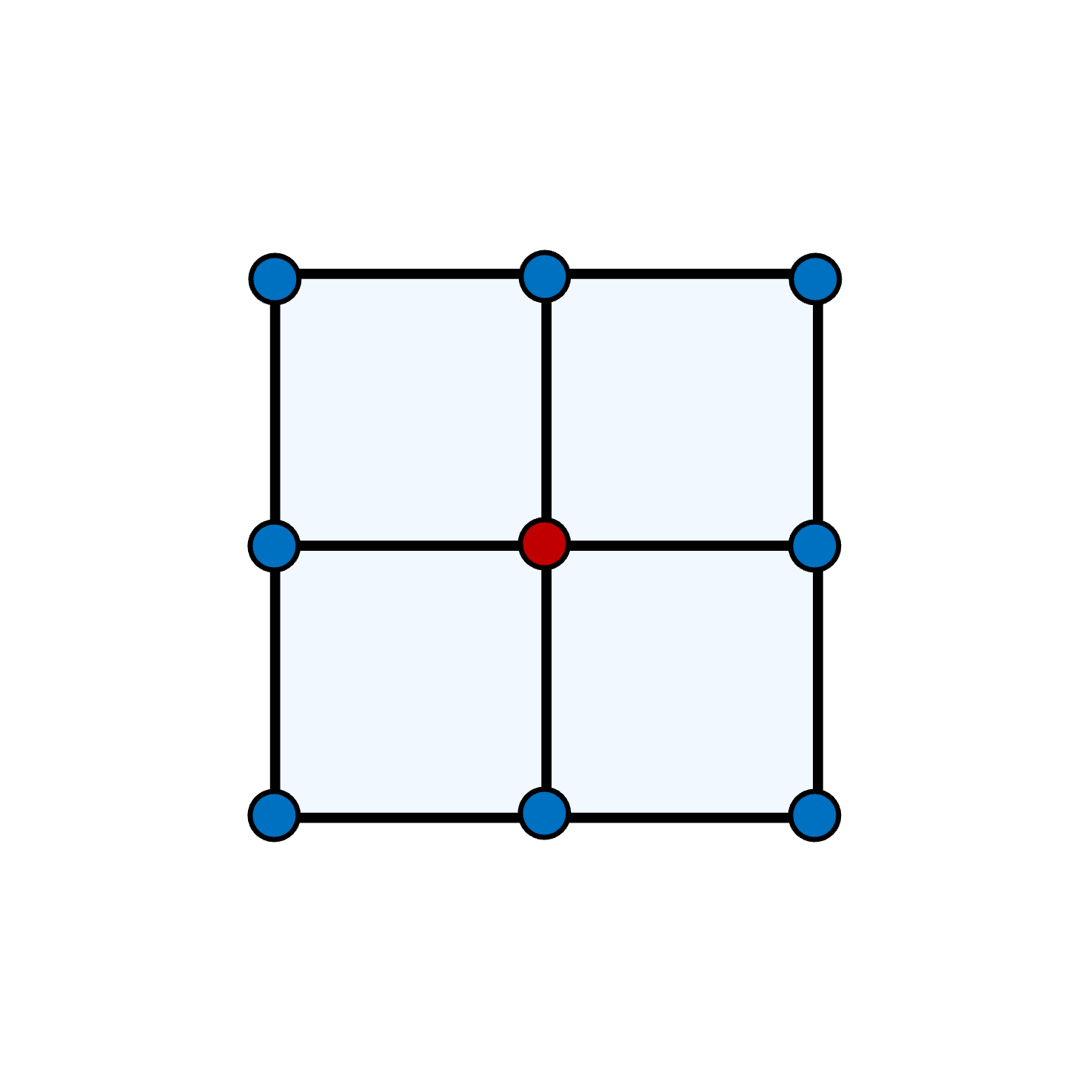} &
\includegraphics[height=0.15\textwidth]{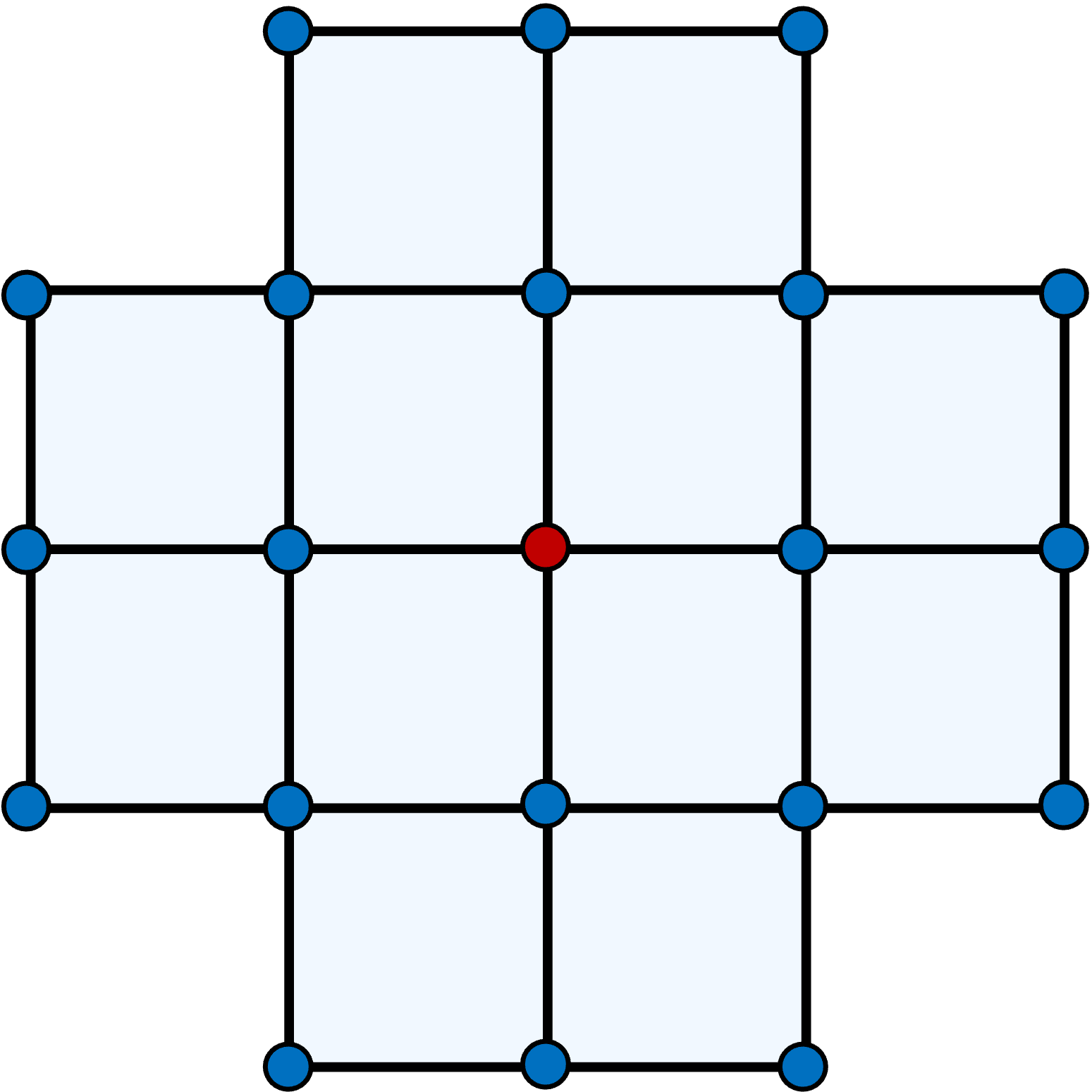} &
\includegraphics[height=0.15\textwidth]{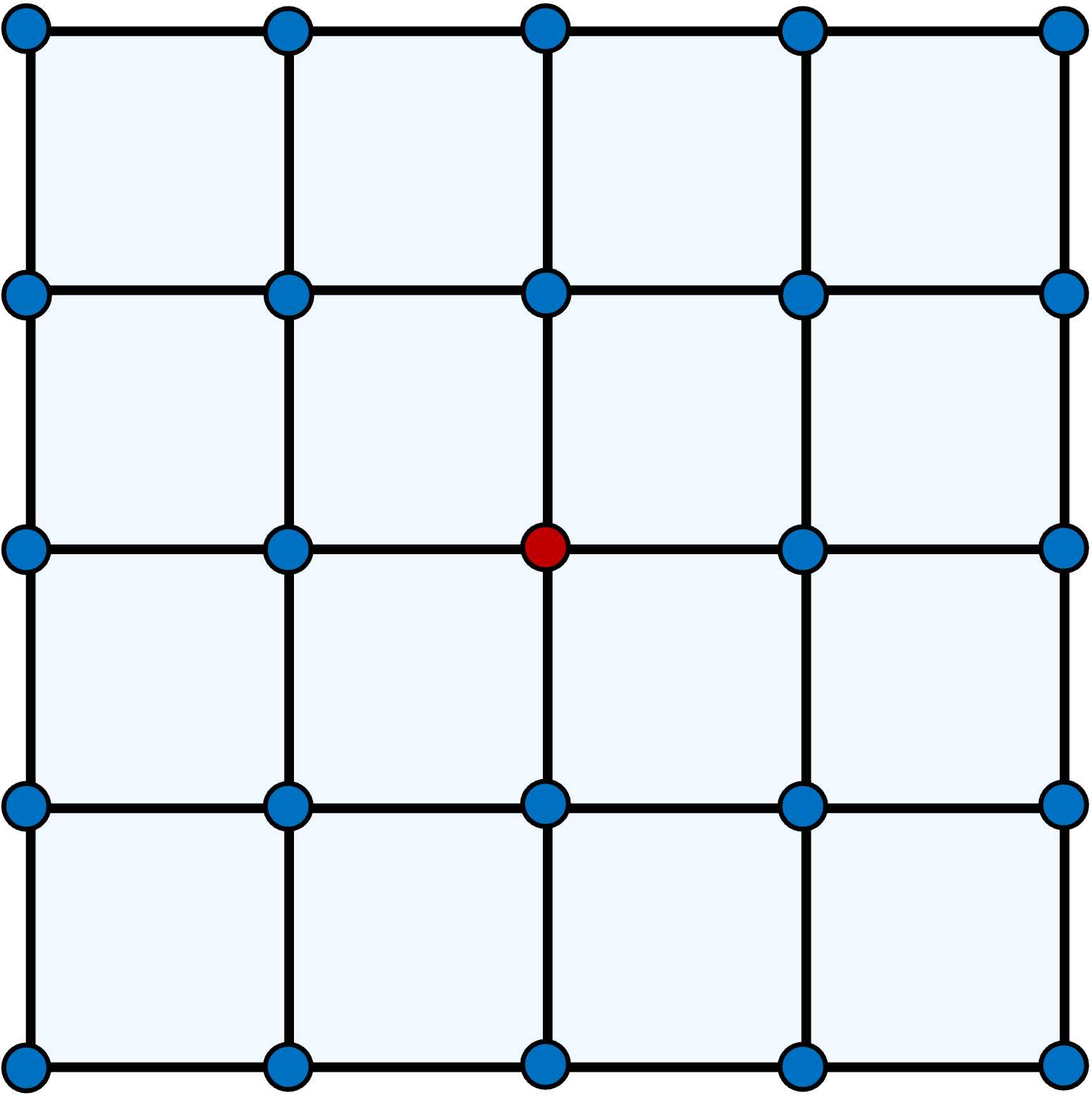}  \\
\scriptsize IFD$(1,1)$ & \scriptsize IFD$(2,2)$ & \scriptsize HIFD$(2,2)$
\end{tabular}
\end{center}
\caption{Stencils used for IFD$(1,1)$, IFD$(2,2)$ and HIFD$(2,2)$ at the red-colored grid.}\label{fig:IFD_stencils}
\end{figure}

We analyze the MILU preconditioners for the IFD and HIFD methods on wide stencils of uniform grids. As illustrated in \cref{fig:IFD_stencils}, IFD and HIFD employ wide stencils, which pose substantial challenges for condition number analysis. 
Considering the 2D Poisson equation \eqref{eq:poisson} on uniform grids with the grid size $h$, implicit finite difference (IFD) \cite{claerbout1985craft} uses a large number of grid points and enhances the accuracy of the finite difference by incorporating a high-order term up to the $(2J+2)$-th order:
\begin{equation*}
\left(1+bh^2\delta_x^2\right)\frac{\partial^2 u}{\partial x^2}= \frac{1}{h^2}\sum_{n=-J}^{J} c_n u\left(x+nh\right) + \cO\left(h^{2J+2}\right),
\end{equation*}
where the adjustable constant $b$ and the coefficients $c_n$ are determined by wave plane theory \cite{liu2009practical}; see detailed derivation in \cite{claerbout1985craft}. Motivated by IFD approach, high-order implicit finite difference (HIFD) adds even higher-order terms, resulting in the following $(2J+4)$-th order accurate formula:
\begin{equation*}
	\left(1+bh^2\delta_x^2 + dh^4\delta_x^4\right)\frac{\partial^2 u}{\partial x^2}= \frac{1}{h^2}\sum_{n=-J}^{J} c_n u\left(x+nh\right) + \cO\left(h^{2J+4}\right).
\end{equation*}
See the detailed derivation in \cite{zapata2017high}. 

For each axis $x$ and $y$, we denote IFD and HIFD that use $J_x$ and $J_y$ points as IFD$\left(J_x,J_y\right)$ and HIFD$\left(J_x,J_y\right)$, respectively. For simplicity, we consider three cases, IFD$(1,1)$, IFD$(2,2)$, and HIFD$(2,2)$, whose coefficient matrix $A$ forms the SPD M-matrix. The values of $b$ and $c_n$ for these cases are summarized in \cref{table:coeff_hifd}. The constant $c_{-n}$ is defined as $c_{-n} = c_n$. See \cref{fig:IFD_stencils} for the stencils for the considered IFD and HIFD schemes.
\begin{table}[htbp]
\footnotesize
	\caption{Coefficients by schemes and orders. }
	\label{table:coeff_hifd}
	\begin{center}
		\begin{tabular}{cc|cccccc}
			\toprule
			& $m$   & $c_0$ & $c_1$ & $c_2$  & $b$ & $d$ \\ \midrule 
			\multirow{2}{*}{IFD}        &  $1$& $-2$  & $1$   &             & $\frac{1}{12}$ &    \\
			&  $2$& $ -\frac{17}{10}$  & $\frac{4}{5}$   &  $\frac{1}{20}$           & $\frac{2}{15}$  &  \\ \midrule 
			\multirow{1}{*}{HIFD}       &  $2$& $ -\frac{53}{42}$  & $\frac{32}{63}$   & $\frac{31}{252}$ &        $\frac{13}{63}$ & $\frac{1}{164}$   \\
			\bottomrule
		\end{tabular}
	\end{center}
\end{table}
Consequently, the linear transform corresponding to the coefficient matrix $A$ is given as follows:
\begin{equation*}\label{eq:high_order_poisson}
	(Au)_{i,j}=   \frac{1}{h^2} \sum_{s=-m}^{m} c_{s} D_y u_{i+s, j}+\frac{1}{h^2} \sum_{l=-m}^{m} c_{s} D_x u_{i, j+s},
\end{equation*}
where $D_x$ and $D_y$ are defined by
\begin{equation*}
	D_xu_{i, s} = \left[
	u_{i-2, s},  u_{i-1, s},  u_{i,s},  u_{i+1, s},  u_{i+2, s}
	\right] v^{\top},
\end{equation*}
and
\begin{equation*}
	D_yu_{s, j} = \left[
	u_{s, j-2},  u_{s, j-1},  u_{s,j},  u_{s,j+1},  u_{s,j+2}
	\right] v^{\top},
\end{equation*}
where the vector $v$ is given by
\begin{equation*}
	v=\left[d, b-4 d, 1-2 b+6 d, b-4 d, d\right].
\end{equation*}
It is easy to verify that $A$ is an SPD M-matrix for the cases under consideration.

The LECN analysis estimates the condition number of the preconditioned system using only local successor and precursor values at each grid. The following Lemma demonstrates how this localized framework significantly simplifies the analysis of the complex structures in IFD and HIFD methods.
\begin{lemma}
Consider $A$ defined using one of methods:~IFD$(1,1)$,~IFD$(2,2)$, and HIFD$(2,2)$. If the value $\tau_{i,j}$ induced by MILU with a lexicographic order, then there exists a constant $c$ such that for any $(i,j)$, the following inequality holds:
    \begin{equation}\label{eq:high_order_tau}
        \tau_{i,j} \le c(4i + j).
    \end{equation}
\end{lemma}
\begin{proof}
    First, note that there exists a constant $c_1$ such that for $(i,j)$ such that $\min(i,j) =1$,
    \begin{equation*}
        \tau_{i,j}\le c_1.
    \end{equation*}
    For $(i,j)$ such that $\min(i,j)\le 3$, by \cref{lemma:tau_evolove}, there exists a constant $c_2$ such that:
    \begin{equation*}
        \tau_{i,j}\le 1 + c_2\max\lrp{\tau_{i-1,j},\tau_{i,j-1}}.
    \end{equation*}
    Thus, there exists a constant $c_3\ge 1$ such that if $\min(i,j)\le 3$, then
    \begin{equation*}
        \tau_{i,j}\le c_3.
    \end{equation*}
    Next, we use mathematical induction on the lexicographic order to prove that for general $i,j\ge4$, $\tau_{i,j}\le c_3(4i+j)$.
    Assume that for $(i',j')$ such that $(i',j')\prec (i,j)$, the induction hypothesis is satisfied.
    Then, by \cref{lemma:tau_evolove}, 
    \begin{equation*}
        \tau_{i,j} \le 1 + \frac{\sum_{K \in s((i,j))} A_{(i,j), K}}{\sum_{K\in p((i,j))} \frac{A_{(i,j), K}}{\tau_{\bi}}}
        \le 1 + \max_{K\in p((i,j))}\tau_{K} \le 1 + c_3(4i+j-1) \le c_3(4i+j).
    \end{equation*}
    Therefore, the induction hypothesis holds.
    This completes the proof.
\end{proof}

The following result directly follows from \cref{proposition:tau_e}.
\begin{corollary}
Consider $A$ defined using one of methods:~IFD$(1,1)$,~IFD$(2,2)$, and~HIFD$(2,2)$. Let $M$ be the MILU preconditioning matrix defined by a lexicographic order. Then, the following equation holds:
\begin{equation*}
    \kappa\lrp{M^{-1}A} = \cO\lrp{h^{-1}}.
\end{equation*}
\end{corollary}

This demonstrates that the LECN analysis can readily assess the condition number of the MILU-preconditioned system for IFD and HIFD methods on wide stencils. This shows the usefulness of the localized approach in the LECN analysis, which relies exclusively on the values of precursors and successors to estimate the condition number, thereby making it highly effective for handling intricate stencil configurations.

\section{Variable Coefficient Poisson Equation on Adaptive Grids}\label{sec:MILU_variable_quad_octree}
We demonstrate the localized approach of the LECN analysis for MILU preconditioner from hierarchical adaptive grids, specifically quadtree and octree. Such a grid can reduce the degrees of freedom required to solve PDEs, saving substantial memory while maintaining comparable accuracy relative to a fully refined grid. In this context, although smaller cells are preferred, the condition number of the coefficient matrix for solving PDEs increases to $\cO\lrp{\bar{h}^{-2}}$, where $\bar{h}$ denote the smallest size of cell. This highlights the importance of designing and analyzing an effective preconditioner on hierarchical adaptive grids to minimize overall computational costs. We also extend the LECN analysis to variable coefficient Poisson equations. For finite volume discretization \cite{frolkovivc2021flux} with Dirichlet boundary conditions on these grids, we prove that the MILU preconditioner achieves a condition number of order $\cO\lrp{\bar{h}^{-1}}$. Experimental validation of the theoretical result is provided in the end. Note that the analyses in this section can be generalized to rectangular (or box) domains with unequal lengths along each axis in 2D (or 3D). However, for clarity, the discussion is foucsed on square (or cubic) domains.


\subsection{Quadtree/Octree}
We define a quadtree in 2D (or octree in 3D) grid and assign a proper order on the related graph to define a MILU preconditioner. The square (or cubic) computational domain is discretized by an uniform grid with square (or cubic) cells as the root of the tree. Then, the doamin is adaptively subdivided into four (or eight) child cells in a quadtree (or octree) where a more accurate numerical solution is necessary.

\begin{definition}[Quadtree/Octree] $\ $

\begin{itemize}
\item[(a)] 
The \textbf{refinement} $\cR(C)$ of a square (or cubic) cell divides it into four equal-sized subcells in 2D (or eight in 3D).

\item[(b)] Let $R_0$ be an unfiorm grid of the computational domain.
A \textbf{quadtree} (or \textbf{octree}) is defined recursively by refining one of the cells $C$ in $R_{i}$ as follows:
For $C\in R_{i}$,
\begin{equation*}
	R_{i+1} = \lrp{R_{i} - C} \cup \cR(C).
\end{equation*}
We refer to $R_0$ as the root grid of $R_{i}$.
\item[(c)] Let the \textbf{length of cell} $C$, $\ell_C$, be the edge length of a square (or cube).
\item[(d)] The \textbf{level} $\mathcal{L}(C) = l$ of cell $C$ in a quadtree/octree grid means that $\ell_C = 2^{-l}h$ where $h$ is the cell length of the root grid.
\item[(e)] 
The \textbf{neighborhood} $n_0\lrp{C}$ of cell $C$ is defined as as the set of all cells that share an edge in 2D (or a face in 3D) with $C$, except for $C$ itself.
\end{itemize}
 \end{definition}
 
\begin{figure}[htbp]
\begin{center}
\begin{tabular}{cc}
\includegraphics[height=0.19\textwidth]{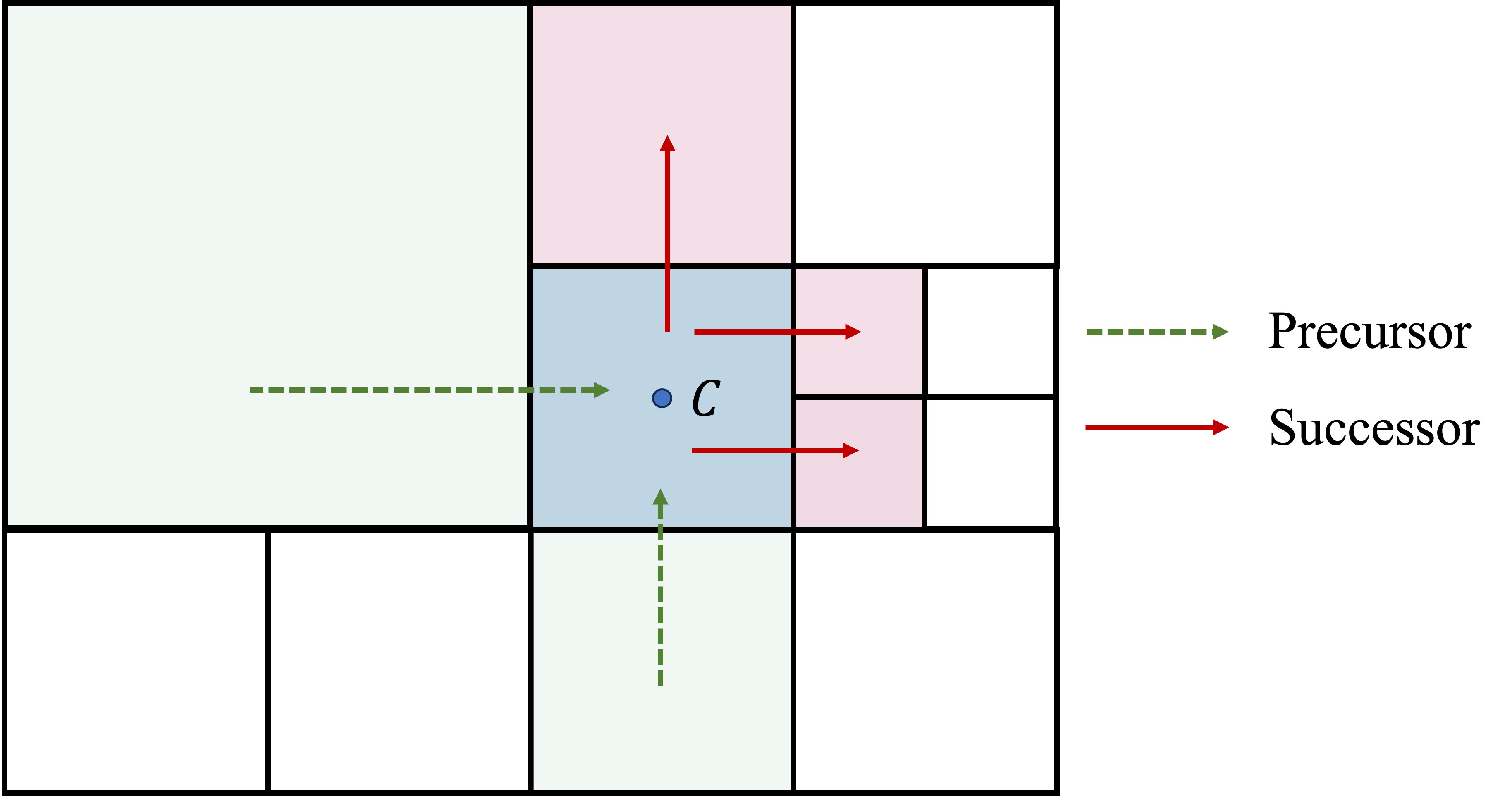} &
\includegraphics[height=0.19\textwidth]{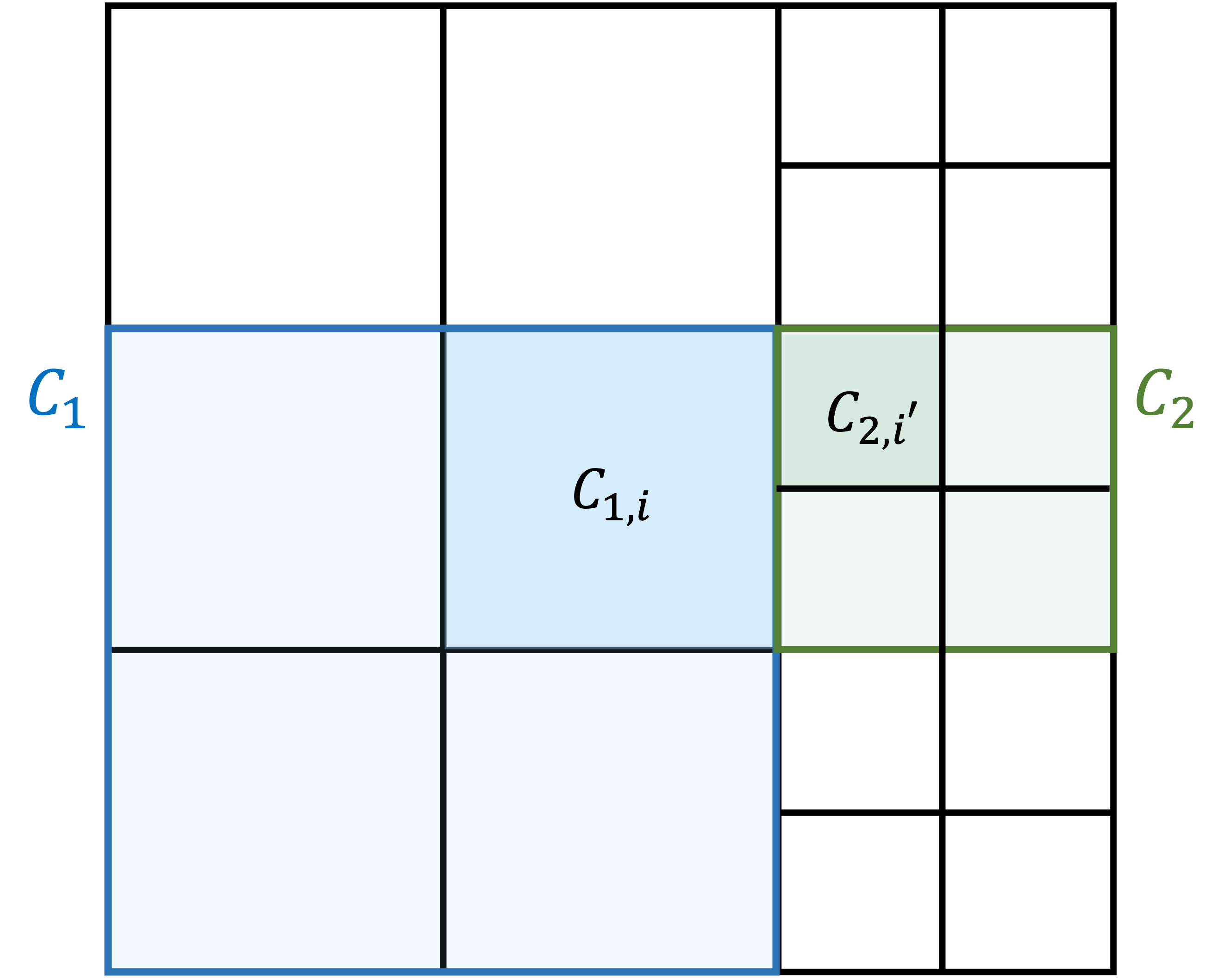} \\
\scriptsize (a) & \scriptsize (b) 
\end{tabular}
\end{center}
 \caption{(a) Illustration of the quadtree order where cells are assigned lexicographical order starting from the bottom-right. For a cell $C$, the two cells highlighted in green are its precursors (a lower order than $C$) and the three neighbor cells highlighted in pink are its successors (a higher order than $C$). (b) An example of a quadtree used to illustrate how the discrepancy in sizes of neighbor cells $C_{1,i}$ and $C_{2,i'}$ causes a jump in the LECN values described in \eqref{eq:example1} and \eqref{eq:example2}.} \label{fig:quadtree_order}
\end{figure}

In the view of graph in~\cref{sec:notation_defintion}, a cell of quadtree (or octree) grid is regarded as a vertex in a graph and neighbor cells are connected by edges in the graph. Now, we define an order on quadtree (or octree) grid, taking into account their hierarchical nature; see~\cref{fig:quadtree_order}-(a) which provides an illustration of the quadtree order. 
\begin{definition}[Order of Quadtree/Octree]\label{def:order_quadtree}
For a quadtree/octree $O$, we define the \textbf{order} $\prec_{O}$ of $O$ as follows:
for a root grid $R_0$ and two cells $C_1, C_2\in R_0$,
\begin{equation*}
    C_1\prec_{R_0} C_2 \iff C_1\prec C_2,
\end{equation*}
where $\prec$ is the lexicographic order for uniform grids.
For $R_{i+1}$ defined as 
\begin{equation*}
    R_{i+1} = \lrp{R_{i} - C} \cup \cR(C),
\end{equation*}
the order between two cells $C_1, C_2\in R_{i+1}$ are defined by,
\begin{equation*}
    C_1\prec_{R_{i+1}} C_2 \iff \begin{cases}
        C_1, C_2\in R_{i} - C \text{ and } C_1 \prec_{R_{i}} C_2, 
        \\ C_1 \in R_{i} - C, C_2\in \cR(C) \text{ and } C_1 \prec_{R_{i}} C, 
        \\ C_1 \in \cR(C) , C_2\in R_{i} - C \text{ and } C \prec_{R_{i}} C_1, 
        \\  C_1, C_2\in \cR(C) \text{ and } C_1 \prec C_2,
    \end{cases}
\end{equation*}
where $\prec$ is the lexicographic order.
For convenience, we will denote all $\prec_{R_i}$ as $\prec_O$. 
\end{definition}

In this study, we limit the size difference between neighbor cells of the quadtree (or octree) to be twice, that is, quadtree (or octree) satisfies the following condition:
\begin{condition}
The level of the two adjacent cells differs by at most one:
if $C_1\in n_0\lrp{C_2}$, 
\begin{equation*}
    \left|\mathcal{L}\lrp{C_1} - \mathcal{L}\lrp{C_2}\right|\le 1.
\end{equation*}
\end{condition}

For simplicity, we introduce notations for the relative sizes of neighbor cells along Cartesian directions. Let $p_{1}$ and $p_{2}$ be the center points of $C_1$ and $C_2$, respectively. For the case $\mathcal{L}\lrp{C_1} = \mathcal{L}\lrp{C_2} +1$, the difference $p \coloneqq p_2 - p_1$ is
\begin{equation*}
    \left(2\ell_{C_1}\right)^{-1} p = \lrp{3,1}, \lrp{3,-1}, \lrp{1,3}, \text{ or } \lrp{-1,3},
\end{equation*}
for the 2D case. In these cases, we will denote 
\begin{align*}
    C_1\prec_{O, x, \pm} C_2 \:\: \text{if} \:\:  \left(2\ell_{C_1}\right)^{-1}p = \lrp{3,\mp 1} \:\: \text{and} \:\: C_1\prec_{O, y, \pm} C_2 \:\: \text{if} \:\:  \left(2\ell_{C_1}\right)^{-1} p = \lrp{\mp 1,3},
\end{align*}
where $x$ and $y$ denote the axis along which $C_1$ and $C_2$ are aligned, and $\pm$ indicates whether the coordinate $p_1$ of the other axis is larger or smaller than that of $p_2$. For the case $\mathcal{L}\lrp{C_1} +1 = \mathcal{L}\lrp{C_2} $, similarly, we denote
\begin{align*}
    C_1\prec_{O, x, \pm} C_2 \:\: \text{if} \:\: \left(2\ell_{C_2}\right)^{-1} p = \lrp{3,\pm 1} \:\: \text{and} \:\: C_1\prec_{O, y, \pm} C_2 \:\: \text{if} \:\: \left(2\ell_{C_2}\right)^{-1} p = \lrp{\pm 1,3}.
\end{align*}
For the 3D case, we define $\prec_{O,\text{axis},\text{sign},\text{sign}}$ in the similar manner. For example, when $\mathcal{L}\lrp{C_1} = \mathcal{L}\lrp{C_2} +1$, $\prec_{O, x, +, +}$ and  $\prec_{O, y, -, +}$ are defined below, respecitvely.
\begin{align*}
    C_1&\prec_{O, x, +, +} C_2 \:\: \text{if} \:\: \left(2\ell_{C_1}\right)^{-1} p = \lrp{3,-1, -1}, \\
    C_1&\prec_{O, y, -, +} C_2 \:\: \text{if} \:\: \left(2\ell_{C_1}\right)^{-1} p = \lrp{1, 3, -1}.
\end{align*}

The goal is to analyze the order of the condition number of the preconditioned system generated by the MILU preconditioner on the mentioned order on a quadtree (or octree). To do so, in the next section, we analyze the asymptotic behavior of the condition number of MILU preconditioner as the cell length is progressively reduced by refined cells into four (resp. eight) cells in a given quadtree (resp. octree) in the following manner:
\begin{definition}[Uniform Refinement of Quadtree/Octree] For a root quadtree (or octree) $O_0$, we recursively define the refinement of $O_i$ as follows:
    \begin{equation*}
        O_{i+1} = \bigcup_{C\in O_i}\mathcal{R}\lrp{C}.
    \end{equation*}
    For each cell $C_k\in O_0$, $C_k$ is divided into $2^n\times 2^n$ or $2^n\times 2^n\times 2^n$ subcells in $O_n$.
For each $\bi\in \{1,\dots, 2^n\}^d$, $C_{k,\bi}$ denotes the $\bi$-th subcell of $C_k$.
\end{definition}

For ease of explanation, we define the index sets $I_0$, $I_1$, and $I_2 \subset \{1,2,\dots, 2^n\}^d$:
\begin{align*}
    I_0 = & \left\{\bi = (i_1,\dots, i_d)\in\{1,2,\dots, 2^n\}^d: 1< i_j < 2^n \text{ for any } j, \right.\\
    &\left. \text{ and } i_j = 2^n-1 \text{ for some } j\in \{1,\dots, d\}   \right\},\\
    I_1 = &\left\{\bi = (i_1,\dots, i_d)\in\{1,2,\dots, 2^n\}^d: i_j= 2^n \text{ for some } j\in \{1,\dots, d\}   \right\},\\
    I_2 = &\left\{\bi = (i_1,\dots, i_d)\in\{1,2,\dots, 2^n\}^d: i_j= 1 \text{ for some } j\in \{1,\dots, d\}   \right\}- I_1.
\end{align*}

\subsection{Finite Volume Method for Variable Coefficient Poisson Equations}\label{sec:fvm2d}
This subsection shows coefficient matrices from the finite volume method (FVM) \cite{frolkovivc2021flux} on quadtree/octree to solve the variable coefficient Poisson equation:
\begin{equation}\label{eq:variable_poisson}
\begin{cases}
   - \nabla\cdot\lrp{\coe \nabla u}  = f \quad &\text{ for } x\in \Omega,
    \\ u = g  \quad &\text{ for } x\in \partial\Omega,
\end{cases}
\end{equation}
where $0<\coe\in C^1\lrp{\Omega, \mathbb{R}}$ is a coefficient function, $f:\Omega\rightarrow \bR$ is a flux function and $g:\partial\Omega\rightarrow\bR$ is a boundary function. A square (or cubic) domain $\Omega$ is discretized by a quadtree (or octree) grid $O = \cup_{i} C_i$, where $C_i$ is a square-shaped (or cube-shaped) cell. A cell $C_i$ of $O$ is a vertex of the graph and $n_0\lrp{C_i}$ is the edge (or face) neighbor of $C_i$. If an edge (or face) of $C_i$ is a subset of $\partial \Omega$, the value $\bar{B}\lrp{C_i}$ indicates the number of edges (or faces) of $C_i$ where the edge (or face) of $C_i$ is a subset of $\partial \Omega$. Then, the coefficient matrix from FVM~\cite{frolkovivc2021flux} to discretize~\eqref{eq:variable_poisson} on the quadtree grid is presented by an SPD M-matrix: 
\begin{equation}\label{eq:fvm2D}
\small
    A_{C_1, C_2} = \begin{cases}
    -2 (\coe_{C_1}^{-1} + \coe_{C_2}^{-1})^{-1} & \text{\small if $C_2\in n_0\lrp{C_1}$, $\mathcal{L}\lrp{C_1} = \mathcal{L}\lrp{C_2}$,}
        \\ - 3\alpha ( \coe_{C_1}^{-1} + 2 \coe_{C_2}^{-1})^{-1} & \text{\small if $C_2\in n_0\lrp{C_1}$, $\mathcal{L}\lrp{C_1} = \mathcal{L}\lrp{C_2} +1$,}        
        \\  -\underset{\tiny C' \in n_0\lrp{C}}{\sum}A_{C,C'}  + 2 \bar{B}(C) \coe_{C} & \text{\small if $C_1 = C_2 (=C)$},
        \\ 0  & \text{\small otherwise},
    \end{cases}
\end{equation}
where $\coe_{C_k}$ is the value of the $\sigma$ at the center of the cell $C_k$ and $\alpha=\frac{2}{3}$ is a constant to indication a common factor of $-A_{C_1,C_2}$ when $\sigma$ is a constant and the condition $|\mathcal{L}\lrp{C_1} - \mathcal{L}\lrp{C_2}| = 1$ holds for neighbor cells $C_1$ and $C_2$. For the octree, the coefficient matrix is obtained similarly:
\begin{equation*}
\small
    A_{C_1, C_2} = \begin{cases} 
    - 2 \ell_{C_1}\lrp{\coe_{C_1}^{-1}+ \coe_{C_2}^{-1}}^{-1} & \text{\small if $C_2\in n_0\lrp{C_1}$, $\mathcal{L}\lrp{C_1} = \mathcal{L}\lrp{C_2}$,}
        \\ - 3\alpha \ell_{C_1} \lrp{\coe_{C_1}^{-1} + 2\coe_{C_2}^{-1}}^{-1} & \text{\small if $C_2\in n_0\lrp{C_1}$, $\mathcal{L}\lrp{C_1} = \mathcal{L}\lrp{C_2} + 1$,}
        \\-\underset{C'\in n_0\lrp{C}}{\sum} A_{C,C'} + 2 \ell_{C} \bar{B}(C)\coe_{C} & \text{\small if $C_1 = C_2 (=C)$},
        \\ 0  & \text{\small otherwise},
    \end{cases}
\end{equation*}
In the following subsection, we apply \cref{proposition:tau_e} to analyze the MILU preconditioner for the coefficient matrices above. Note that since $\coe$ is a $C^1$-function, there exists a constant $P_{\coe}$ and an integer $n$ such that for any $C_1\in \mathcal{V}$ and $C_2\in n_0\lrp{C_1}$,
\begin{equation}\label{eq:variable_bound}
    \left|\frac{\coe_{C_2}}{\coe_{C_1}}-1\right|\le   \frac{P_{\coe}}{2^{n+1}}.
\end{equation}

\subsection{LECN Analysis on Adaptive Grids}\label{sec:fvm2d_MILU}
In this subsection, we apply LECN analysis to examine the condition number of the MILU preconditioner for the FVM method introduced in the previous subsection. The main results, presented in \cref{theorem:2d_octree_main}, \cref{corollary:main}, and \cref{theorem:3d_octree_main}, show the condition number of preconditioned system from MILU preconditioner on quadtree or octree grid is improved to the order $\cO\left(\bar{h}^{-1}\right)$ for the smallest cell length $\bar{h}$. 

\begin{theorem}\label{theorem:2d_octree_main}
    For a quadtree $O_0$ defined on the square domain $[0, a_x]\times [0,a_y]$ and $a_x = a_y$, let 
    \begin{equation*}
       O_0 = \left\{C_1, \dots, C_m\right\}.
    \end{equation*}
    Assume that each cell $C_k\in O_0$ is defined as 
    \begin{equation*}
        C_k = [a_{k,x} - \ell_{C_k}, a_{k,x}] \times [a_{k,y} - \ell_{C_k}, a_{k,y}].
    \end{equation*}
    Then, the following relation holds for any $n\in \mathbb{N}$, $\bi\in I_0$ and $C_{k,\bi}\in O_n$,
    \begin{equation}\label{eq:main_I1}
        \tau_{k,\bi} \le g_{k,n, \bi} + o\lrp{2^n}, 
    \end{equation}
    where
    \begin{equation*}
        g_{k,n,\bi} = \exp\lrp{\lrp{a_{k,x}+a_{k,y}- \frac{\lrp{2^{n+1} -i - j}h}{2^{\mathcal{L}\lrp{C_k} + n}}}\sup\frac{\|\nabla \coe\|_1}{\coe}}\frac{a_{k,x}+a_{k,y}}{h}  2^{n + \mathcal{L}\lrp{C_k}+1},
    \end{equation*}
    and there exists a constant $c\in \mathbb{R}$ independent of $n$, that satisfies the following inequality for any $n\in \mathbb{N}$ and $\bi$:
    \begin{equation*}
        \tau_{k,\bi} \le cg_{k,n,\bi} + o\lrp{2^n}.
    \end{equation*}
\end{theorem}
In conjunction with \cref{proposition:tau_e}, we directly obtain the following corollary.
\begin{corollary}\label{corollary:main}
    For a quadtree $O_0$ and $n$-times refinement $O_n$, consider the MILU preconditioning matrix $M$. 
    Then, the following equation holds:
    \begin{equation*}
        \kappa\lrp{M^{-1}A} = \cO\lrp{2^n}.
    \end{equation*}
    Furthermore, if the coefficient function $\coe$ is constant, there exists a constant $c$ independent of $n$ and the quadtree $O_0$ such that
        \begin{equation*}
        \kappa\lrp{M^{-1}A} \le  c 2^{n + \max_{C\in O_0}\mathcal{L}\lrp{C}}.
    \end{equation*}
\end{corollary}

\cref{theorem:2d_octree_main} and \cref{corollary:main} imply that the condition number increases with the linear growth as the depth of the quadtree increases. We also check the same growth numerically in the experiments in \cref{subsec:quadtree_exp}.

We emphasize that that the current results are generalization of the previous results in \cite{yoon2015analyses} in the view of variable coefficients and hierarchical adaptive grids. The challenge of generalization is mainly caused by the non-uniformity in the number and length of neighbor cells. It complicates the analysis of the condition number, estimating LECN $\tau_{k,\bi}$ on quadtree (or octree), when a boundary of cell is shared with multiple smaller cells due to the presence of T-junction. To illustrate this, let us consider a simple case where $\sigma=1$. In~\cref{fig:quadtree_order}-(b), we have two neighbor cells $C_1$ and $C_2$ with $\mathcal{L}\lrp{C_1}+1 = \mathcal{L}\lrp{C_{2}}$ and then cells $C_{1,\bi}$ and $C_{2,\bi'}$ are the child cells of $C_1$ and $C_2$, respectively,  after $C_1$ and $C_2$ have been once refined. In this case, the LECN values $\tau$ with $\bi=\left(i_1,i_2\right)$ and $\bi'=\left(i'_1,i'_2\right)$ for the two cells $C_{1,\bi}$ and $C_{2,\bi'}$ are provided as follows:
\begin{align}
        \tau_{1,\bi} &= 1 + \frac{1 + 2\alpha}{\frac{1}{\tau_{1,\lrp{i_1-1,i_2}}} + \frac{1}{\tau_{1,\lrp{i_1,i_2-1}}}}, \label{eq:example1} \\
        \tau_{2,\bi'} &= 1+ \frac{2}{\frac{1}{\tau_{2,\lrp{i'_1,i'_2-1}}} + \frac{\alpha}{\tau_{1,\bi}}}. \label{eq:example2}
\end{align}
This indicates that when neighboring cells have different lengths, a jump in $\tau$ occurs.
{If the bounds on $\tau$ are estimated without rigorous consideration, one might erroneously estimate suboptimal amplification factors of  $\frac{1+2\alpha}{2}$ and $\frac{2}{1+\alpha}$ for $\tau_{1,\bi}$ and $\tau_{2,\bi'}$, respectively, which cause the unbounded growth of $\tau$. (or potentially indicating unbounded behavior of $\tau$.) However, thorough quantitative analysis demonstrates that the actual amplification factors are $2\alpha$ and $\frac{1}{\alpha}$, respectively, when cells of different lengths are neighbor each other. This analysis eventually confirms the bounded nature of $\tau$ under rigorously defined conditions.} Accordingly, the proof of \cref{theorem:2d_octree_main} needs Lemmas \ref{lemma:jump_double}, \ref{lemma:jump_half}, and \ref{lemma:2d_path_finding} and is presented after~\cref{lemma:2d_path_finding}. The following Lemmas \ref{lemma:jump_double} and \ref{lemma:jump_half} provide bounds for $\tau$ at the subcells located on the edges of the pre-refined cells where the levels of neighbor cells increase and decrease, respectively. The proofs of two lemmas are provided in \cref{subsec:proof:jump_double} and \cref{subsec:proof:jump_half}, respectively.
\begin{lemma}\label{lemma:jump_double}
Consider a quadtree $O$ and a cell $C_k$.
Assume there exists two cells $C_{k_1}, C_{k_2}\in p\lrp{C_k}$ such that
\begin{equation*}
    \mathcal{L}\lrp{C_{k_1}} = \mathcal{L}\lrp{C_{k_2}} = \mathcal{L}\lrp{C_{k}} + 1, \  C_{k_1}\prec_{O, x, -} C_k \text{ and }  C_{k_2}\prec_{O, x, +} C_k.
\end{equation*}
Define $\eta_{i}$ as
\begin{equation*}
    \eta_{i} = \begin{cases}
        \tau_{k_{1}, \lrp{2^n-1, i}} &\text{ if } 1\le i\le 2^n,
       \\ \tau_{k_{2}, \lrp{2^n-1,i-2^n}} &\text{ if } 2^n+1\le i\le 2^{n+1}.
    \end{cases}
\end{equation*}
   Assume there exists a constant $b, c\in \mathbb{R}_+$ such that 
\begin{equation*}
   0\le  \eta_{i} \le c2^n + b,
\end{equation*}
for any $ 1\le i\le 2^{n+1} -1$ and $n\in \mathbb{N}$.
If $2^{n-2}\ge P_{\coe}$, then there exist constants $c_1, c_2$ satisfying the followings: for any $n\in \mathbb{N}$ and $1\le i\le 2^n-1$,
\begin{align*}
\begin{split}
  \tau_{k,\lrp{1,i}}\le\lrp{1+\frac{2P_{\coe}}{2^n}}\lrp{1 + \frac{P_{\coe}\lrp{\alpha+1}}{\alpha 2^n}}\frac{c2^n+b}{2} + \lrp{c_12^n + c_2} \lrp{\frac{3}{4}}^{-i} +  9.
   \end{split}
\end{align*}
\end{lemma}

\begin{lemma}\label{lemma:jump_half}
Consider a quadtree $O$ and a cell $C_k$.
Assume that there exists two cells $C_{k_1}, C_{k_2}\in s\lrp{C_k}$ such that
\begin{equation*}
    \mathcal{L}\lrp{C_{k_1}} = \mathcal{L}\lrp{C_{k_2}} = \mathcal{L}\lrp{C_{k}} + 1,\  C_k\prec_{O, x, -}  C_{k_1} \text{ and }  C_k\prec_{O, x, +}  C_{k_2}.
\end{equation*}
Assume there exists constants $b, c\in \mathbb{R}_+$ such that 
\begin{equation*}
   0\le \tau_{k,\bi} \le c2^n + b,
\end{equation*}
for any $\bi = \lrp{2^n-1,i}$ where $1\le i\le 2^n-1$ and $n\in \mathbb{N}$.
Define $\beta_{i}$ as
\begin{equation*}
    \beta_{i} = \begin{cases}
        \tau_{k_{1}, \lrp{1, i}} &\text{ if } 1\le i\le 2^n,
       \\ \tau_{k_{2}, \lrp{1,i-2^n}} &\text{ if } 2^n+1\le i\le 2^{n+1}.
    \end{cases}
\end{equation*}
Then, for $2^{n-2}\ge P_{\coe}$, there exist constants $c_1, c_2$ such that 
\begin{equation*}
   0\le \beta_{i} \le 2\lrp{1+\frac{2P_{\coe}}{N}} \lrp{1 + \frac{\lrp{1+2\alpha}P_{\coe}}{2\alpha N}}(c2^n+b) + 5+ \lrp{c_12^n+c_2} \lrp{\frac{4}{5}}^{-i},
\end{equation*}
for any $n\in \mathbb{N}$ and $1\le i\le 2^{n+1}-1$.
\end{lemma}

A technical lemma is presented below for handling the recurrence relations such as those in \eqref{eq:example1} and \eqref{eq:example2}, whose proof is provided in \cref{subsec:proof:tau_recurrent_beta}.
\begin{lemma}\label{lemma:tau_recurrent_beta}
Consider constants $\alpha,b_1,c_1,b_2, c_2, b_3,c_3,\gamma\in \mathbb{R}_+$ such that $\gamma< 1$ and $N\in \mathbb{N}$.
        Assume that the following inequality holds for $\beta_1$: 
    \begin{equation*}
     \beta_{1}\le c_1N + b_1.
    \end{equation*}
For $i\ge 2$, assume the following recurrence relation: 
    \begin{equation*}
        \beta_{i} \le 1 + \frac{1+\alpha}{\frac{1}{\beta_{i-1}} + \frac{1}{c_2N+b_2 + \lrp{c_3N+b_3}\gamma^{i}}}.
    \end{equation*}
    Then, for any $\gamma_2$ satisfying
    \begin{equation*}
      1>  \gamma_2 \ge \gamma \text{ and } \gamma_2 > \frac{1}{1 +\alpha},
    \end{equation*}
      there exist constants $\alpha_1, \alpha_2\in \mathbb{R}_+$ such that for any $i\in \mathbb{N}$ and $N\in \mathbb{N}$, the following inequality holds:
    \begin{equation}\label{eq:2d_inequality_beta}
\beta_{i} 
     \le \alpha \lrp{c_2 N  + b_2} + \lrp{\alpha_1 N + \alpha_2}\gamma_2^{i}+\frac{2+\alpha}{\alpha}.
    \end{equation}
    Furthermore, if there exists a constant $u$ such that $\alpha\le u$, then, $\alpha_1$ and $\alpha_2$ do not depend on $\alpha$.
\end{lemma}

The previous lemmas establish that if the $\tau$ values of subcells within a cell $C_k$ are bounded by a specific value, then $\tau_{k', (1,i_2)}$ decays exponentially, converging to double or half the value as $i_2$ increases. However, this could potentially lead to large $\tau$ values in the remaining regions. Fortunately, combinatorial techniques ensure that the influence of large $\tau_{k',(1,i_2)}$, $\tau_{k',(i_1,1)}$ values for small $i_1$ and $i_2$ on $\tau_{2^n-1,i_2}$ and $\tau_{i_1, 2^n-1}$   diminishes due to the exponential convergence. 
The following \cref{lemma:2d_path_finding} addresses this phenomenon and the proof is provided in \cref{subsec:proof:2d_path_finding}.
\begin{lemma}\label{lemma:2d_path_finding}
Consider constants $c, c_1, c_2, \gamma, b\in \mathbb{R}_+$ such that $\gamma<1$.
    Assume that for $\bi=\lrp{i,j}$, if $i=1$ and $2\le j<2^n$, then,
    \begin{equation*}
        \tau_{k,\bi} \le c2^n + \lrp{c_12^n+c_2}\gamma^j +b.
    \end{equation*}
    Similarly, assume that if $j=1$ and $2\le i<2^n$, then,
    \begin{equation*}
        \tau_{k,\bi} \le c2^n + \lrp{c_12^n+c_2}\gamma^i +b.
    \end{equation*}
    Then, for $\bi = \lrp{i,j}\in I_0$, the following inequality holds:
    \begin{equation*}
        \tau_{k, \lrp{i,j}} \le \exp\lrp{ {\frac{h(i+j)}{2^{\mathcal{L}\lrp{C_k} + n}} \sup\frac{\|\nabla \coe \|_1}{\coe} }}\lrp{c2^n + i + j} +o\lrp{2^n}.
    \end{equation*}
\end{lemma}

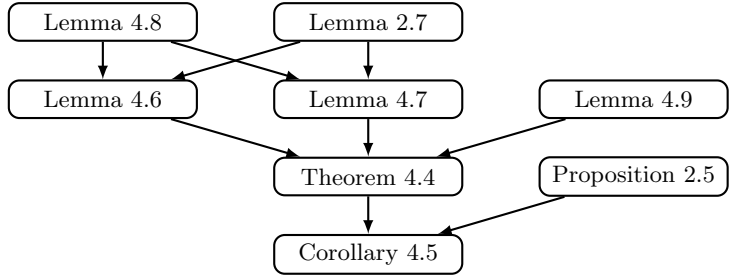
\begin{figure}[htbp]
    \centering    
\begin{tikzpicture}[
    font=\small,
    thick,
    cus_block/.style={
        draw,
        rectangle,
        rounded corners,
        minimum width=2.5cm,
        minimum height=0.5cm,
        align=center}]

\node[cus_block] (block1) {\cref{lemma:tau_recurrent_beta}};
\node[cus_block, right=1cm of block1] (block2) { \cref{lemma:tau_evolove}};
\node[cus_block, below=0.5cm of block1] (block3) { \cref{lemma:jump_double}};
\node[cus_block, below=0.5cm of block2] (block4) { \cref{lemma:jump_half}};
\node[cus_block, right=1cm of block4] (block5) { \cref{lemma:2d_path_finding}};
\node[cus_block, below=0.5cm of block4] (block6) { \cref{theorem:2d_octree_main}};
\node[cus_block, right=1cm of block6] (block7) { \cref{proposition:tau_e}};
\node[cus_block, below=0.5cm of block6] (block8) { \cref{corollary:main}};
\draw[-latex] (block1) edge (block3)
(block1) edge (block4)
    (block2) edge (block3)
    (block2) edge (block4)
    (block3) edge (block6)
    (block4) edge (block6)
    (block5) edge (block6)
    (block6) edge (block8)
    (block7) edge (block8);

\end{tikzpicture}
    \caption{A logical flow in the proof of \cref{theorem:2d_octree_main} and \cref{corollary:main}, illustrating the progression from Lemmas and \cref{proposition:tau_e} to the main results on condition number estimation for MILU preconditioners on quadtree.}\label{fig:pf_flow}
\end{figure}
By combining these results, we can conclude the proof of the theorem; please see \cref{fig:pf_flow} for an overview of the overall logical flow of the proof.

\begin{proof}[Proof of \cref{theorem:2d_octree_main}]
    Using mathematical induction on the quadtree order of $O_0$, we will prove Eq (\ref{eq:main_I1}) for any $\bi\in I_0$.
    First, consider $C_k\in O_0$ such that 
    \begin{equation*}
        C_{k'} \nprec_0 C_k , 
    \end{equation*}
    for any $C_{k'}$. 
    The left and lower boundary of subcells $C_k$ are also the boundaries of the entire quadtree, and thus $\tau_{k,(1,j)}, \tau_{k,(i,1)} = o(2^n)$. 
    By \cref{lemma:2d_path_finding}, for $\bi = (i_1, i_2)\in I_0$, 
    \begin{align*}
    \begin{split}
        \tau_{k,\bi} &\le 
        \exp\lrp{ {\frac{h(i+j)}{2^{\mathcal{L}\lrp{C_k}+n}} \sup\frac{\|\nabla \coe\|_1}{\coe} }}(i+j) +o\lrp{2^n} \le g_{k,n, \bi}+ o\lrp{2^n}.
        \end{split}
    \end{align*}
    Now consider a general $k$ and assume that Eq \cref{eq:main_I1} holds for any $k'$ such that $C_{k'}\prec_O C_k$ and $\bi\in I_0$.
    For the vertical edge of $C_k$ adjacent to its precursor, there 
are three cases:
    \begin{itemize}
        \item $\mathcal{L}\lrp{C_{k_1}} = \mathcal{L}\lrp{C_{k_2}} = \mathcal{L}\lrp{C_{k}}+1$ and ${C_{k_1}}\prec_{O, x, +} C_k$ and $\mathcal{L}\lrp{C_{k_2}}\prec_{O, x, -} C_k$.
        By the induction hypothesis and \cref{lemma:jump_double}, there exist constants $c_1, c_2$ such that for any $i\le 2^n-1$,
         \begin{align*}
  \begin{split}
         \tau_{k,\lrp{1,i}} \le & \frac{g_{k_1, n, (2^n, 2^n)}}{2} + \lrp{c_12^n+c_2}\lrp{\frac{3}{4}}^{i} +  o\lrp{2^n} .
        \end{split}
        \end{align*}
           \item $\mathcal{L}\lrp{C_{k_1}} = \mathcal{L}\lrp{C_{k}}$.
                   By the induction hypothesis, with straightforward calculations, there exist constants $c_1, c_2$ such that the following inequality holds for any $i\le 2^n-1$:
         \begin{align*}
  \begin{split}
           \tau_{k,\lrp{1,i}}
            \le & g_{k_1,n, (2^n, 2^n)} + \lrp{c_12^n+ c_2}\lrp{\frac{4}{5}}^{i} +  o\lrp{2^n}.
        \end{split}
        \end{align*}
    \item $\mathcal{L}\lrp{C_{k_1}} +1 = \mathcal{L}\lrp{C_{k}}$.
    If ${C_{k_1}}\prec_{O, a, -} C_k$, by the induction hypothesis and \cref{lemma:jump_half}, there exist constants $c_1,c_2$ such that for any $i\le 2^n-1$,
            \begin{align*}
  \begin{split}
           \tau_{k,\lrp{1,i}}\le 2g_{k_1,n, (2^n, 2^n)} + \lrp{c_12^n+ c_2}\lrp{\frac{4}{5}}^{i} +  o\lrp{2^n}.
       \end{split}
       \end{align*}         
       If ${C_{k_1}}\prec_{O, a, +} C_k$, similarly,
       \begin{align*}
           \begin{split}
                \tau_{k,\lrp{1,i}}\le 2g_{k_1,n, (2^n, 2^{n-1})} + \lrp{c_12^n+ c_2}\lrp{\frac{4}{5}}^{i} +  o\lrp{2^n}.
           \end{split}
       \end{align*}
    \end{itemize}
   Define $\beta $ as $2^{n-1}$ if ${C_{k_1}} +1 = \mathcal{L}\lrp{C_{k}}$ and $\mathcal{L}\lrp{C_{k_1}}\prec_{O, a, +} C_k$, and as zero otherwise.  
   Then, the following inequality holds for any $\bi=(i,j)\in I_0$:
         \begin{align*}
  \begin{split}
       &\exp\lrp{ {\frac{(i+j)h}{2^{\mathcal{L}\lrp{C_k} +n}} \sup\frac{\|\nabla \coe\|_1}{\coe} }}\lrp{2^{\mathcal{L}(C_k) - \mathcal{L}(C_{k_1})}g_{k_1,n,(2^n, 2^n - \beta)} + i+j} 
       \\ &\le  \exp\lrp{ {\frac{(i+j)h}{2^{\mathcal{L}\lrp{C_k} +n}} \sup\frac{\|\nabla \coe\|_1}{\coe} }}\exp\lrp{\lrp{a_{k_1,x}+a_{k_1,y}- \frac{\beta h}{2^{\mathcal{L}\lrp{C_k} + n}}}\sup\frac{\|\nabla \coe\|_1}{\coe}}
       \\& \times \lrp{2^{\mathcal{L}(C_k) - \mathcal{L}(C_{k_1})} 
      \frac{a_{k_1,x}+a_{k_1,y}}{h}  2^{n + \mathcal{L}\lrp{C_{k_1}}+1} + i+j} 
      \\ &\le  \exp\lrp{\lrp{a_{k,x}+a_{k,y}- \frac{\lrp{2^{n+1} -i - j}h}{2^{\mathcal{L}\lrp{C_k} + n}}}\sup\frac{\|\nabla \coe\|_1}{\coe}}
      \lrp{      \frac{a_{k,x}+a_{k,y}}{h}  2^{n + \mathcal{L}\lrp{C_{k}}+1}}
      \\ &\le g_{k,n,\bi}.
   \end{split}
   \end{align*}
    Similar results hold for $i\le 2^n-1$ and $\tau_{k,\lrp{i,1}}$.
    By \cref{lemma:2d_path_finding}, for $\bi=(i,j)\in I_0$,
     \begin{align*}
  \begin{split}
      \tau_{k,\lrp{i,j}}       \le g_{k,n,\bi} +o\lrp{2^n}.
   \end{split}
   \end{align*}
    Thus, the induction hypothesis is satisfied for any $C_k$ and Eq \cref{eq:main_I1} holds for $\bi\in I_0$.
    We can easily verify that
    \begin{align*}
        \tau_{k,\bi}\le \lrp{1+2\alpha}\max_{\bi'\in I_0}\tau_{k,\bi'} + o\lrp{2^n}, & \text{ for } \bi\in I_1,\\
        \tau_{k,\bi}\le 2\max_{k'\in p\lrp{k}, \bi'\in I_1}\tau_{k',\bi'} + o\lrp{2^n}, & \text{ for } \bi\in I_2.
    \end{align*}
    For general $\bi=\lrp{i,j}$ such that $1\le i,j\le 2^n-1$,
    \begin{equation*}
        \tau_{k,\bi}\le\exp\lrp{ {\frac{2h}{2^{\mathcal{L}\lrp{C_k}}} \sup\frac{\|\nabla \coe\|_1}{\coe} }}\max_{\bi'\in I_2}\lrp{\tau_{k,\bi'} +i+j} + o\lrp{2^n}.
    \end{equation*}
    Combining these results, there exists a constant $c$ such that the following inequality holds:
    \begin{equation*}
        \tau_{k,\bi} \le c\exp\lrp{\lrp{a_{k,x}+a_{k,y}}\sup\frac{\|\nabla \coe\|_1}{\coe}}\frac{a_{k,x}+a_{k,y}}{h}  2^{n + \mathcal{L}\lrp{C_k}+1} + o\lrp{2^n}.
    \end{equation*}
    This completes the proof.
\end{proof}

Similar to the quadtree case, we perform LECN analysis of the MILU preconditioner on an octree in 3D. The main theorem for octrees is stated below and we provide the proof in \cref{subsec:proof:3d_octree_main}. 
Note that, unlike its 2D counterpart, the theorem does not provide a bound for the constant $c$ in terms of the octree structure $O_0$ or the coefficient function $\coe$. 
However, it does provide that $\tau$ is bounded by $\cO\lrp{2^n}$.
\begin{theorem}\label{theorem:3d_octree_main}
    For an octree $O_0$ in the domain $[0, a_x]\times [0,a_y]\times [0,a_z]$ defined as $       O_0 = \left\{C_1, \dots, C_m\right\}$,
    there exists a constant $c$ such that satisfy the following inequality for any $n\in \mathbb{N}$ and $\bi$:
    \begin{equation*}
        \tau_{k,\bi} \le c \cdot 2^n.
    \end{equation*}
\end{theorem}

\subsection{Numerical Validation}\label{subsec:quadtree_exp} We provide a numerical validation to support the theoretical analyses discussed in \cref{sec:fvm2d_MILU} regarding the improvement of the condition number of the matrices \eqref{eq:fvm2D}. The numerical results also highlight an influence on the computational efficiency. The implementation in C++ is executed on a computer with a 3.40 GHz CPU and 32.0 GB of RAM. The power method is used to compute the condition number of the matrix, while all linear systems are solved by the preconditioned conjugate gradient (PCG) method. 

\begin{figure}[htbp]
\begin{center}
\begin{tabular}{cc}
\includegraphics[width=0.25\textwidth]{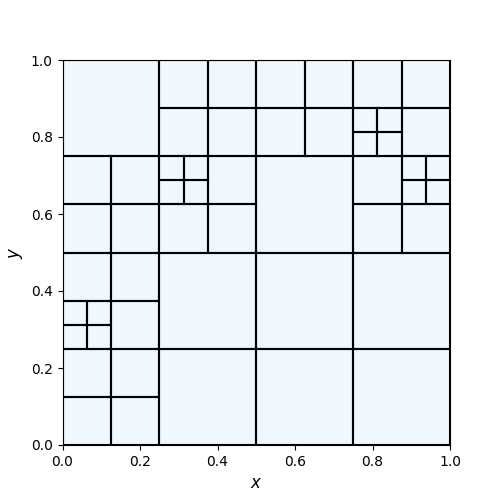} &
\includegraphics[width=0.25\textwidth]{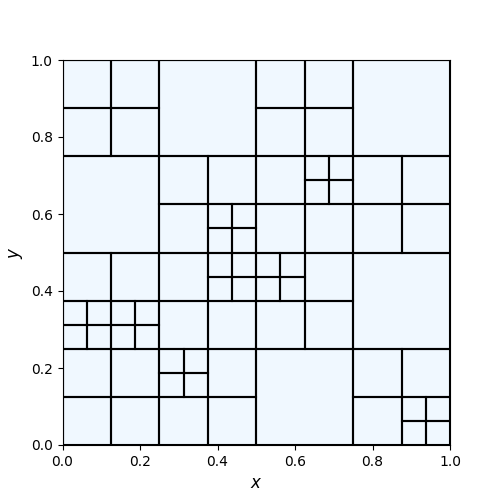} \\
\scriptsize (a) Example 1 & \scriptsize (b) Example 2
\end{tabular}
\end{center}
\caption{Visualization of randomly generated initial quadtrees.}
\label{fig:rand_quadtrees}
\end{figure}

\paragraph{Condition Number Reduction}
Two elliptic coefficients $\sigma$ in \eqref{eq:variable_poisson} are tested:
\begin{itemize}
    \item Example 1: An oscillatory coefficient $\sigma_1\left(x,y\right)=\sin\left(\pi x\right)\cos\left(2\pi y\right) + 1.5$.
    \item Example 2: A steep coefficient $\sigma_2\left(x,y\right)=\exp\left(3-2x\right)y\left(3-3y\right)+0.5$.
\end{itemize}
To validate the MILU while ruling out any possible coincidence arising from the regularity of the grid, we randomly generate quadtree grids depicted in \cref{fig:rand_quadtrees}, which are then subsequently refined. For comparison, the Jacobi and ILU preconditioners are tested on the same examples on quadtree. 



\begin{figure}[htbp]
\begin{center}
\begin{tabular}{cc}
\includegraphics[width=0.30\textwidth]{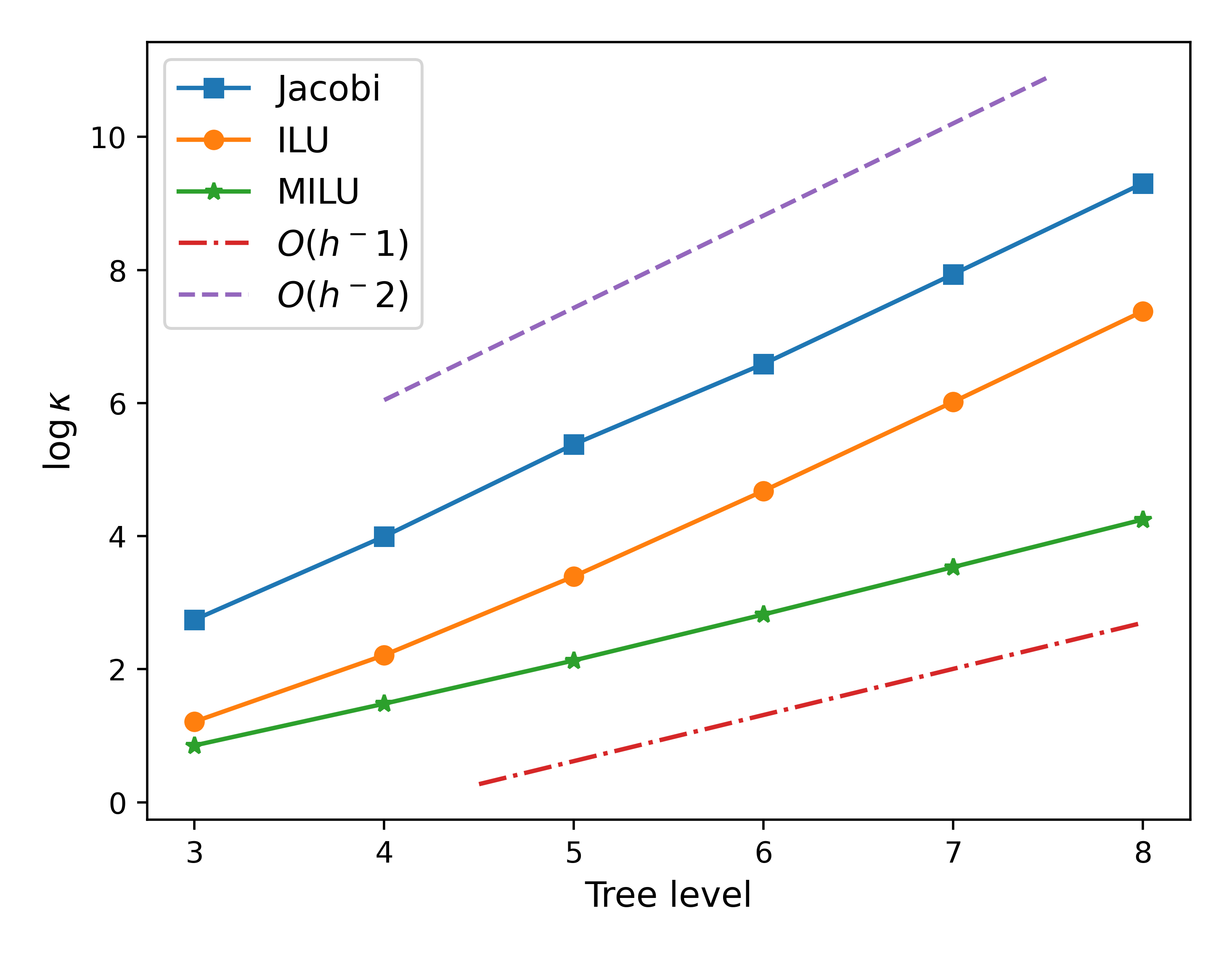} &
\includegraphics[width=0.30\textwidth]{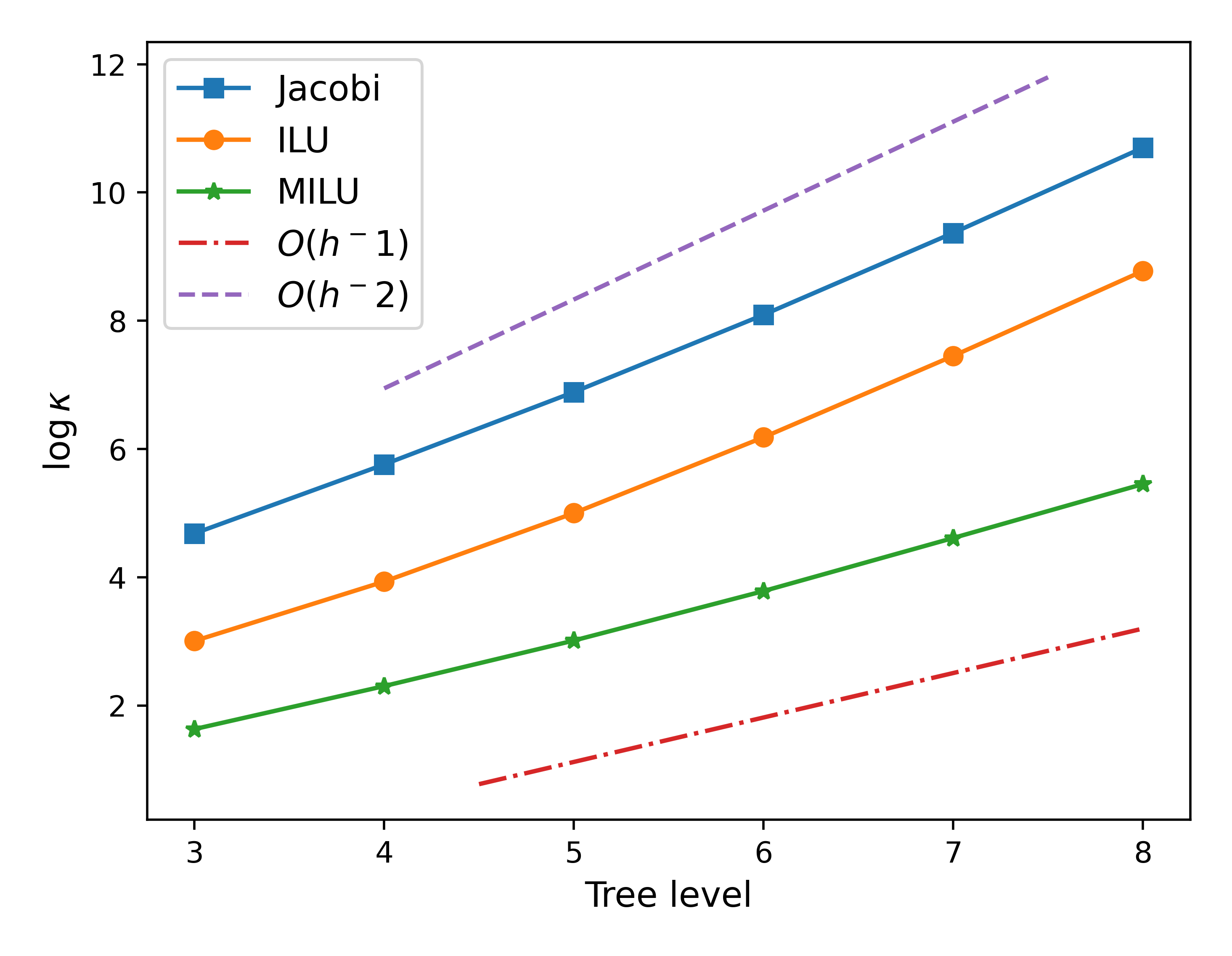} \\
\scriptsize (a) Example 1 & \scriptsize (b) Example 2
\end{tabular}
\end{center}
\caption{Comparison of $\kappa\left(M^{-1}A\right)$ between Jacobi, ILU, and MILU preconditioners in the random quadtree. The MILU preconditioner reduces the condition number to $\cO\left(\bar{h}^{-1}\right)$.}
\label{fig:quadtree}
\end{figure}

In \cref{fig:quadtree}, the numerical results are shown by log-scale plots of the tree level versus the condition number. The results align with the theoretical analysis in \cref{theorem:2d_octree_main}. The MILU preconditioner reduces the condition number by order $\cO\left(\bar{h}^{-1}\right)$ in both examples, outperforming the Jacobi and ILU preconditioners along the condition number of order $\cO\left(\bar{h}^{-2}\right)$.


\begin{figure}[htbp]
\begin{center}
\begin{tabular}{c}
\includegraphics[width=0.40\textwidth]{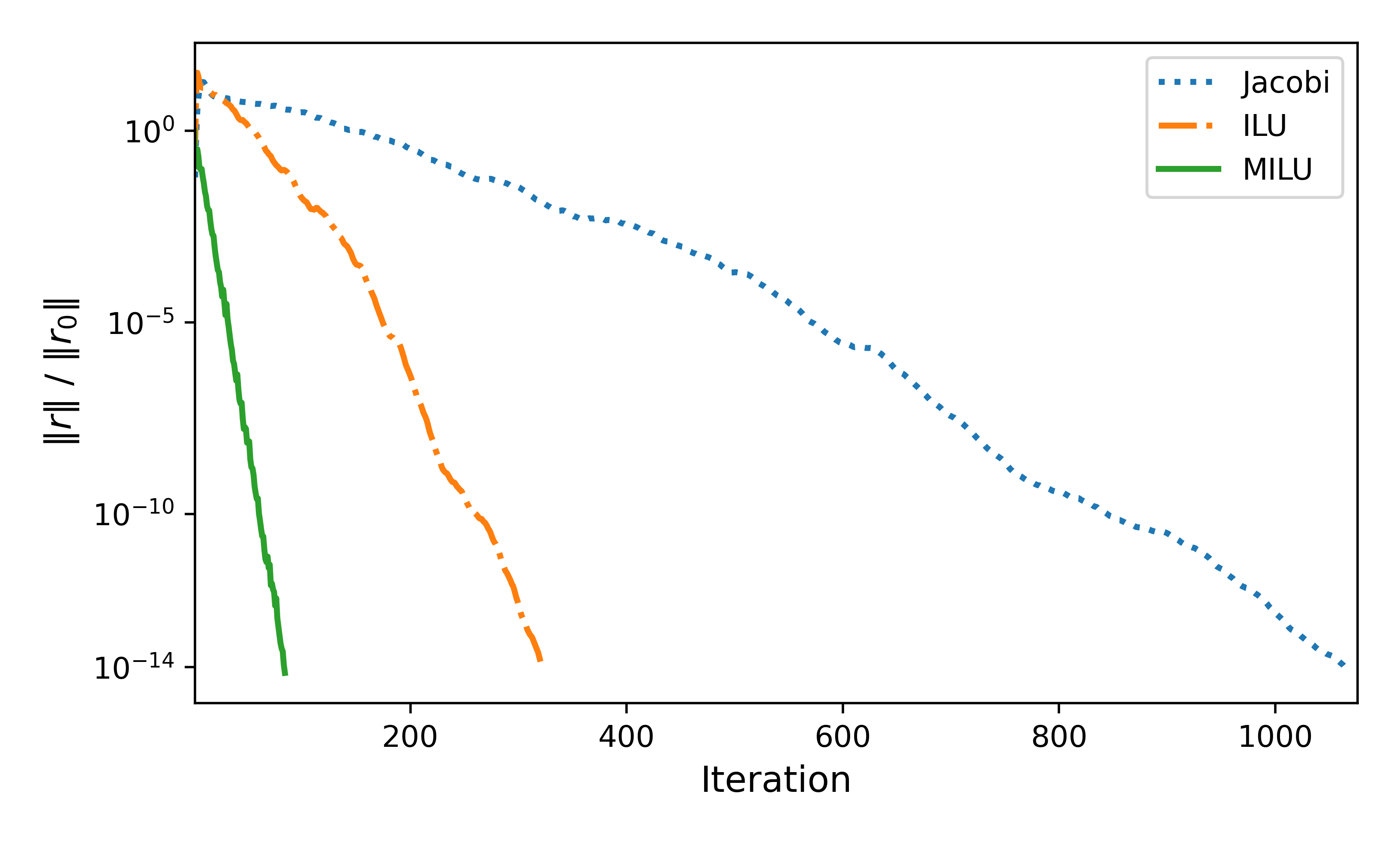}
\end{tabular}
\end{center}
\caption{The plot shows the number of iterations versus the relative residual. The results indicate that a preconditioner with a smaller condition number speeds up convergence, as seen in the fewer iterations needed to achieve the same relative residual error.}
\label{fig:pcg}
\end{figure}

\paragraph{Practical Impact on PCG Convergence}
In the view of the practical advantages of MILU in solving linear systems, the number of iterations in the preconditioned conjugate gradient (PCG) method is checked to reach a relative residual error $10^{-14}$. It is tested by Example 2 on the domain discretized with a randomly generated quadtree grid with depth $8$ with $74752$ cells. In \cref{fig:pcg}, a notable improvement of the number of iterations in the PCG convergence is presetned when the MILU preconditioner is used. Specifically, MILU requires only about $8\%$ as many iterations as Jacobi and $26\%$ as many as ILU, demonstrating its significantly higher efficiency in terms of iteration count and computational cost.

\section{Conclusion}
In this paper, we introduced a framework for analyzing Modified Incomplete LU (MILU) preconditioners, expanding their applicability by defining them a weighted undirected graph with self-loops. A key contribution is the \textit{Localized Estimator of Condition Number (LECN)}, a novel measure that facilitates an estimation of condition number by local measurements, thereby simplifying the analysis of complex and previously unexplored matrix structures. The usability of the LECN analysis was demonstrated through three cases: revisiting estabilished results on uniform grids, analyzing high-order finite difference schemes on wide stencils, and extending the framework to hierarchical adaptive grids such as quadtree and octree for variable coefficient Poisson equations. This work demonstrates the potential of LECN analysis to address challenging problems in numerical linear algebra and PDE solvers, providing a powerful tool for estimating condition number in many applications.

A possible direction of future research is the expansion of the LECN framework to a broader class of PDEs, unstructred grids such as polyhedral meshes, and Neumann boundary conditions. From a theoretical standpoint, a deeper analysis to establish a tight upper bound in the LECN would be a valuable investigation. On the practical side, exploring the optimal order on graphs to enhance the efficiency of the MILU preconditioner is crucial, as our findings demonstrate that an order of graph provides insights into estimating the condition number. Finally, while extending LECN analysis to coefficient matrices that are not M-matrices presents current challenges, overcoming these could enable the application of the framework to a broader spectrum of matrix structures.

\appendix
\section{Technical Proofs}
\subsection{Proof of \cref{lemma:tau_recurrent_beta}}\label{subsec:proof:tau_recurrent_beta}
\begin{proof}
As $1 + \frac{1+\alpha}{\frac{1}{x} + \frac{1}{c_2N+b_2 + \lrp{c_3N+b_3}\gamma^{i}}}$ is an increasing function with respect to $x$, we only need to consider the case $\beta_1\ge \alpha\lrp{c_2N+b_2}$.
  First, prove that $\beta_{i}\ge \alpha\lrp{c_2N+b_2}$ for any $i\in \mathbb{N}$.
  We use mathematical induction.
    When $i=1$, there is nothing to prove.
    For $i_0\in \mathbb{N}$, assume the induction hypothesis is satisfied for any $i< i_0$.
    Then,
     \begin{align*}
  \begin{split}
        \beta_{i_0} &= 1 + \frac{1+\alpha}{\frac{1}{\beta_{i_0-1}}+ \frac{1}{c_2N+b_2 + \lrp{c_3N+b_3}\gamma^{i-1}}}
        \ge \frac{1+\alpha}{\frac{1}{\beta_{i_0-1}}+ \frac{1}{c_2N+b_2}}
       \\ 
       &\ge\frac{1+\alpha}{\frac{1}{\alpha\lrp{c_2N+b_2}}  + \frac{1}{c_2N+b_2}}
        = \alpha\lrp{c_2N+b_2} .
    \end{split}
    \end{align*}
    Thus, the induction hypothesis holds for any $i\in \mathbb{N}$, and therefore, $\beta_i \geq \alpha\lrp{c_2N+b_2}$ for all $i \in \mathbb{N}$.

 Now we prove the upper bound.
     As $\gamma_2 \ge \gamma$ and $\gamma_2 > \frac{1}{1+\alpha}$, we can select constants $1>\epsilon_1, \epsilon_2\in \mathbb{R}_+$ such that
     \begin{equation*}
        \frac{1}{1 +\alpha} <   \frac{1}{1 +\frac{\alpha}{1+\epsilon_1}} < \gamma_2, \text{ and } 
          1 + \alpha\epsilon_2 < \gamma_2\lrp{ {1 +\frac{\alpha}{1+\epsilon_1}} }.
     \end{equation*}
     Then, the following inequality holds for any $\alpha'\ge \alpha$:
     \begin{equation*}
           \frac{1 +\epsilon_2\alpha'}{1 + \frac{\alpha'}{1+\epsilon_1}}\le \gamma_2.
     \end{equation*}
     Consider a natural number $m\in \mathbb{N}$ such that if $k>m$, the following inequality holds:
          \begin{equation*}
         \max\lrp{\frac{c_3}{c_2},\frac{b_3}{b_2}} {\gamma}^m\le \epsilon_1.
      \end{equation*}
    Note that $\epsilon_1, \epsilon_2$, and $m$ do not depend on $N$.
  Define $\alpha_1, \alpha_2$, and $b$ as \begin{align*}
      \begin{split}
          \alpha_1& = \max\lrp{\frac{c_1}{\gamma_2}, \lrp{1+\alpha}\lrp{\frac{c_2}{\gamma_2^{m}} +  c_3}, \frac{\alpha c_3}{\epsilon_2}
  },\\
  \alpha_2 & = \max\lrp{\frac{b_1}{\gamma_2}, \lrp{1 +\alpha}\lrp{\frac{b_2 +1}{\gamma_2^{m}} + b_3}, \alpha b_2 +\frac{2+\alpha}{\alpha}, \frac{\alpha b_3}{\epsilon_2}},\\
  b & = \alpha b_2 + \frac{2+\alpha}{\alpha}.
      \end{split}
  \end{align*}
    Define $q_{k}$ as 
    \begin{equation*}
      q_{k} =   \alpha{c_2 N} + \lrp{\alpha_1 N + \alpha_2}\gamma_2^{k} + b.
    \end{equation*}
We want to prove that for any $i\in \mathbb{N}$,
\begin{equation*} 
    \beta_{i} \le q_{i}.
\end{equation*}
 Use mathematical induction on $i$ to prove that Eq (\ref{eq:2d_inequality_beta}) holds.
        If $i=1$, as $c_1\le \gamma_2\alpha_1$ and $b_1\le \gamma_2\alpha_2$,
    \begin{equation*}
        \beta_1\le c_1N + b_1 \le    \lrp{\alpha c_2  + \gamma_2\alpha_1}N +\lrp{ b + \gamma_2\alpha_2} = q_1. 
    \end{equation*}
    Therefore, for $i=1$, the induction hypothesis is satisfied.
    For $1<i\le m$, as $ \lrp{1+\alpha}\lrp{\frac{c_2}{\gamma_2^{m}} + c_3}\le \alpha_1$ and $ \lrp{\alpha + 1}\lrp{\frac{b_2+1}{\gamma_2^{m}} + b_3}\le \alpha_2$,
    \begin{equation*}
        \beta_{i}\le 1 + \lrp{1+\alpha}\lrp{c_2N+b_2 + \lrp{c_3N+b_3}\gamma^i} \le \alpha c_1N+b+ \lrp{\alpha_1 N+\alpha_2}\gamma_2^{i} = q_{i}. 
    \end{equation*}
    Thus, the induction hypothesis is satisfied for $i\le m$.
    
    For $i > m$, assume the induction hypothesis is satisfied for $i-1$.
    By the induction hypothesis,
    \begin{equation*}
        \beta_{i} 
         \le 1 +\frac{1+\alpha}{\frac{1}{\beta_{i-1}} + \frac{1}{c_2N+b_2 + \lrp{c_3N+b_3}\gamma^{i}}} 
      \le   1 +\frac{1+\alpha}{\frac{1}{q_{i-1}}  + \frac{1}{c_2N+b_2 + \lrp{c_3N+b_3}\gamma^{i}}} 
    .
    \end{equation*}
    Thus,
     \begin{align*}
  \begin{split}
        & \beta_{i} - \alpha c_2N - b\\
        &
        \le 1 +\frac{1+\alpha}{\frac{1}{q_{i-1}} + \frac{1}{c_2N+b_2 + \lrp{c_3N+b_3}\gamma^{i}}}  -\alpha c_2N - \alpha b_2 +\lrp{\alpha b_2 - b}
      \\  &\le 1 + \lrp{\alpha b_2 -b} + \frac{\lrp{1 -\frac{\alpha c_2N + b }{q_{i-1}}} +\alpha \lrp{1 -\frac{c_2N+b_2}{c_2N+b_2 + \lrp{c_3N+b_3}\gamma^{i}}} 
      + \lrp{b -\alpha b_2}\lrp{\frac{1}{q_{i-1}}}
      }{\frac{1}{q_{i-1}} + \frac{1}{c_2N+b_2 + \lrp{c_3N+b_3}\gamma^{i}}} 
       \\  &\le 1 + \lrp{\alpha b_2 -b} + \frac{\lrp{1 -\frac{\alpha c_2N + b }{q_{i-1}}} +\alpha \lrp{1 -\frac{c_2N+b_2}{c_2N+b_2 + \lrp{c_3N+b_3}\gamma^{i}}} 
      + \lrp{b -\alpha b_2}\lrp{\frac{1}{q_{i-1}}} 
      }{\frac{1}{q_{i-1}} + \frac{1}{(1+\epsilon_1) \lrp{c_2N+b_2}}} 
    \\ & \le 1 + \lrp{\alpha b_2 -b} + \frac{ \frac{\lrp{\alpha_1N+\alpha_2}\gamma_2^{i-1}}{q_{i-1}}   + \frac{\alpha\lrp{c_3N+b_3}\gamma^i}{c_2N+b_2} }{\frac{1}{q_{i-1}} + \frac{1}{(1+\epsilon_1) \lrp{c_2N+b_2}}}
    + \frac{\lrp{b-\alpha  b_2}\lrp{\frac{1}{q_{i-1} }}}{\frac{1}{q_{i-1}}  +\frac{1}{ (1+\epsilon_1) \lrp{c_2N+b_2}}}
    \\ & \le 1  - \lrp{b -\alpha b_2}\lrp{1 - \frac{1}{1 +\frac{\alpha}{2}}} 
    + \lrp{\alpha_1N+\alpha_2}\gamma_2^{i-1}
    \lrp{ 
    \frac{1 +\epsilon_2\frac{q_{i-1}}{c_2N+b_2}}{1 + \frac{q_{i-1}}{(1+\epsilon_1) \lrp{c_2N+b_2}}}
    }
    \\ & \le \lrp{\alpha_1N+\alpha_2}\gamma_2^{i}.
    \end{split}
    \end{align*}
    Therefore, 
    \begin{equation*}
        \beta_{i}\le \alpha c_2N + b + \lrp{\alpha_1N+\alpha_2}\gamma^{i} = q_i,
    \end{equation*}
    and the induction hypothesis is satisfied.
    We conclude that $\beta_{i}\le q_{i}$, and this completes the proof.
\end{proof}

\subsection{Proof of \cref{lemma:jump_double}}\label{subsec:proof:jump_double}
\begin{proof}
Let $N = 2^n$. By \cref{lemma:tau_evolove} and Eq (\ref{eq:variable_bound}),
\begin{align*}
    \tau_{k_1, \lrp{2^n, 1}} 
    &\le   1 +\frac{\sum_{K_2\in s\lrp{C_{k_1, \lrp{2^n, 1}}}}-A_{\lrp{C_{k_1, \lrp{2^n, 1}}},K_2}}{\sum_{K_1\in p\lrp{C_{k_1, \lrp{2^n, 1}}}} \frac{-A_{\lrp{C_{k_1, \lrp{2^n, 1}}},K_1}}{\tau_{K_1}}} \\
    &
    \le 1 + \lrp{1+ \frac{P_{\coe}}{N}}\lrp{2\alpha+1}{\tau_{k_1, \lrp{2^n-1, 1}}}
    \\ &\le \lrp{2\alpha+1}c\lrp{1+ \frac{P_{\coe}}{N}} N + \lrp{2\alpha+1}b\lrp{1+ \frac{P_{\coe}}{N}}+1.
\end{align*}
For $i=1,\dots, 2^{n+1}$ define $\alpha_i$ as 
\begin{equation*}
    \alpha_i = \begin{cases}
        \tau_{k_1, \lrp{2^n, i}} &\text{ if } 1\le i\le 2^n,
       \\ \tau_{k_2, \lrp{2^n, i-2^n}} &\text{ if } 2^n+1\le i\le 2^{n+1}.
    \end{cases}
\end{equation*}
Then, $\alpha_1 \le \lrp{1+ \frac{P_{\coe}}{N}}\lrp{2\alpha+1}\lrp{cN+b} + 1$.
By \cref{lemma:tau_evolove} and (\ref{eq:variable_bound}), we have the following recurrent inequality:
\begin{equation*}
    \alpha_{i+1}\le 1 + \lrp{1 + \frac{P_{\coe}}{N}}\frac{\alpha+1}{\frac{1}{\alpha_i} + \frac{1}{cN+b}} = 1 + \frac{\alpha + \frac{P_{\coe}\lrp{\alpha+1}}{N}+1}{\frac{1}{\alpha_i}+ \frac{1}{cN+b}} .
\end{equation*}
By \cref{lemma:tau_recurrent_beta}, there exist constants $c_1, c_2$, and $1> \gamma> \frac{1}{\alpha + \frac{P_{\coe}\lrp{\alpha+1}}{N}+1}$ such that
 \begin{align*}
  \begin{split}
    \alpha_i &\le \lrp{\alpha + \frac{P_{\coe}\lrp{\alpha+1}}{N}} (cN + b) + \frac{\alpha+2}{\alpha} + \gamma^{i} \lrp{c_1 N+ c_2}\eqqcolon \alpha'_i.
    \end{split}
    \end{align*}
Define $\beta_i$ as
\begin{equation*}
    \beta_{i} = \tau_{k, \lrp{1, i}}.
\end{equation*}
Then, by \cref{lemma:tau_evolove} and Eq (\ref{eq:variable_bound}),
\begin{equation*}
    \beta_1 \le 1 + \frac{2\lrp{1+\frac{P_{\coe}}{N}}}{\frac{\alpha}{\alpha_1} + \frac{\alpha}{\alpha_2}} \le 1 + \frac{\lrp{1+\frac{P_{\coe}}{N}}\alpha_1}{\alpha}\le  1 + \lrp{1+\frac{P_{\coe}}{N}}^2\lrp{2\alpha+1}\frac{cN + b + 1}{\alpha}.
\end{equation*}
 Again, by \cref{lemma:tau_evolove} and Eq (\ref{eq:variable_bound}), we have the following recurrent inequality:
 \begin{equation*}
     \beta_{i+1}\le 1+ \frac{2\lrp{1+\frac{P_{\coe}}{N}}}{\frac{1}{\beta_i} + \frac{\alpha}{\alpha_{2i+1}} + \frac{\alpha}{\alpha_{2i+2}}}
     \le 1+ \frac{2\lrp{1+\frac{P_{\coe}}{N}}}{\frac{1}{\beta_i} + \frac{2\alpha}{\alpha'_{2i+2}}}
     .
 \end{equation*}
 By \cref{lemma:tau_recurrent_beta}, there exist constants $c'_1, c'_2$ independent of $i$ and $n$ such that the following inequality holds:
\begin{align*}
    \beta_{i} \le & \lrp{1+\frac{2P_{\coe}}{N}}\lrp{1 + \frac{P_{\coe}\lrp{\alpha+1}}{\alpha N}}\frac{cN + b}{2} \\ &+ \frac{3}{2}\lrp{\frac{\alpha+2}{\alpha} } + \frac{1+2}{1} +  \lrp{c_1' + c_2'N}\lrp{\frac{3}{4}}^{i}
    \\ 
    \le & \lrp{1+\frac{2P_{\coe}}{N}}\lrp{1 + \frac{P_{\coe}\lrp{\alpha+1}}{\alpha N}}\frac{cN+b}{2} +  \lrp{c_1' + c_2'N}\lrp{\frac{3}{4}}^{i} +  9.
\end{align*}
This completes the proof.
\end{proof}

\subsection{Proof of \cref{lemma:jump_half}}\label{subsec:proof:jump_half}

\begin{proof}
Let $N = 2^n$.
    For $i=1,\dots,2^n$, define $\alpha_i$ as 
    \begin{equation*}
        \alpha_i = \tau_{k,\lrp{2^n,i}}.
    \end{equation*}
    Then, by \cref{lemma:tau_evolove} and Eq \eqref{eq:variable_bound},
    \begin{equation*}
        \alpha_1 \le 1 + \lrp{1 +\frac{P_{\coe}}{N}}\lrp{cN + b}\lrp{1 + 2\alpha},
    \end{equation*}
    and
    \begin{equation*}
        \alpha_{i+1}\le 1+ \lrp{1 +\frac{P_{\coe}}{N}}\frac{1 + 2\alpha}{\frac{1}{\alpha_{i}} + \frac{1}{cN+b}}.
    \end{equation*}
    By \cref{lemma:tau_recurrent_beta}, there exist constants $c_1, c_2$, and $\gamma > \frac{1}{1 + 2\alpha}$ such that
    \begin{equation*}
        \alpha_i\le \lrp{2\alpha + \frac{\lrp{1+2\alpha}P_{\coe}}{N}} (cN + b)+ \frac{2\alpha+2}{2\alpha} +\lrp{c_1N+c_2}\gamma^i.
    \end{equation*}
    For $i =1,\dots,2^{n+1}$, define $\beta_i$ as 
    \begin{equation*}
        \beta_i = \begin{cases}
        \tau_{k_1, \lrp{1, i}} &\text{ if } 1\le i\le 2^n,
       \\ \tau_{k_2, \lrp{1, i-2^n}} &\text{ if } 2^n+1\le i\le 2^{n+1}.
       \end{cases}
    \end{equation*}
    Then, by \cref{lemma:tau_evolove} and Eq \eqref{eq:variable_bound},
    \begin{equation*}
        \beta_1 \le 1+ \lrp{1+\frac{P_{\coe}}{N}}\frac{2 }{\frac{\alpha}{\alpha_1}}\le 1 + \lrp{1+\frac{P_{\coe}}{N}} + \lrp{1+\frac{P_{\coe}}{N}}^2\frac{\lrp{2cN + 2b}\lrp{1+2\alpha}}{\alpha},
    \end{equation*}
    and
    \begin{equation*}
        \beta_{i+1}\le 1+\frac{2\lrp{1+\frac{P_{\coe}}{N}}}{\frac{1}{\beta_i} + \frac{\alpha}{\alpha_{\floor{\frac{i+1}{2}}}}}.
    \end{equation*}
     By \cref{lemma:tau_recurrent_beta}, there exist constants $c'_1, c'_2$, and $\gamma > \frac{1}{1 + 2\alpha}$ such that
         \begin{align*}
      \beta_i &\le  2\lrp{1+\frac{2P_{\coe}}{N}} \lrp{1 + \frac{\lrp{1+2\alpha}P_{\coe}}{2\alpha N}}(cN + b) + 5 + \lrp{c'_1N+c'_2}\lrp{\frac{4}{5}}^{i}
    .
    \end{align*}
    This completes the proof.
\end{proof}

\subsection{Proof of \cref{lemma:2d_path_finding}}\label{subsec:proof:2d_path_finding}
\begin{proof}
Let $N = 2^n$.
First, prove the case when $\coe$ is a constant function: $\coe(x)=1$ for any $x\in \mathbb{R}^2$.
 Consider a constant $m = \ceil{-\log_{\gamma}\lrp{\lrp{c_1N+c_2}/\lrp{1-\gamma}}}$ so that if $l\ge m$, we have
 \begin{equation*}
     \lrp{c_1N+c_2}\gamma^l < 1 - \gamma.
 \end{equation*}
 Note that $m =\cO\lrp{n}$.
 If $l\ge m$, the following inequality holds:
 \begin{equation*}
     \tau_{k,\lrp{1,l}}\le cN + b+1.
 \end{equation*}
 If $l < m$,  then
 \begin{equation*}
      \tau_{k,\lrp{1,l}}\le cN + b + \lrp{c_1N + c_2}.
 \end{equation*}
     For $2\le i,j\le N-1$, the following inequality holds:
 \begin{equation*}
     \tau_{k,\lrp{i,j}}\le 1 + \frac{2}{\frac{1}{\tau_{k,\lrp{i-1,j}}} +\frac{1}{\tau_{k,\lrp{i,j-1}}} }\le 1+ \frac{\tau_{k,\lrp{i-1,j}} + \tau_{k,\lrp{i,j-1}}}{2}.
 \end{equation*}
 Thus,
 \begin{equation*}
     \tau_{k,\lrp{i,j}} - i - j \le \frac{\lrp{\tau_{k,\lrp{i-1,j}} - i -j+1} + \lrp{\tau_{k,\lrp{i,j-1}} - i -j +1}}{2}.
 \end{equation*}
Define $\psi_{i,j}$ as 
 \begin{equation*}
     \psi_{i,j} = \tau_{k,\lrp{i,j}} - i - j .
 \end{equation*}
 Then,
 \begin{equation*}
      \psi_{i,j}\le \frac{\psi_{i-1,j}+ \psi_{i,j-1}}{2}.
 \end{equation*}
 We will prove that 
 \begin{equation*}
     \psi_{i,j} \le \sum_{l=2}^i   \psi_{l,1}  \xi\lrp{i,j,l} + \sum_{l=2}^j   \psi_{1,l}  \xi\lrp{i,j,l},
 \end{equation*}
 where $\xi\lrp{i,j,l}\coloneqq 2^{-i-j+l+1} \binom{i+j-l- 2}{j-2}$.
 We use mathematical induction on the lexicographic order of $\bi =\lrp{i,j}$.
For $\lrp{i,j} = \lrp{2,2}$,
\begin{equation*}
    \psi_{2,2} \le \frac{ \psi_{1,2} +  \psi_{2,1} }{2},
\end{equation*}
and the induction hypothesis is satisfied.
 Assume that for $j<j_0$ and $\bi = \lrp{2, j}$, the induction hypothesis is satisfied:
 \begin{equation*}
     \psi_{2,j_0-1}\le   \psi_{2,1}   2^{-j_0+2}
     + \sum_{l=2}^{j_0-1}   \psi_{1,l}   2^{-j_0 +l}.
 \end{equation*}
 Then,
  \begin{align*}
  \begin{split}
     \psi_{2,j_0} &\le \frac{ \psi_{1,j_0} +  \psi_{2,j_0-1}  }{2}
\le \frac{\psi_{1,j_0}}{2} + \frac{ \psi_{2,1}  2^{-j_0+2}
     + \sum_{l=2}^{j_0-1}   \psi_{1,l}  2^{-j_0 +l}}{2}\\
&     = \psi_{2,1}  2^{-j_0+1} + \sum_{l=2}^{j_0}\psi_{1,l}  2^{-j_0 -1+l}.
     \end{split}
 \end{align*}
Thus, the induction hypothesis is satisfied for all $(2, j)$.
Similar result holds for $(i,2)$.
Then, for any $\lrp{i,j}$ such that $i_0,j_0\ge 3$, assume that $\lrp{i_0-1,j_0}$ and $\lrp{i_0,j_0-1}$ satisfy the induction hypothesis.
For $i,j\ge 3$, 
 \begin{align*}
  \begin{split}
    \psi_{i_0,j_0}  &\le   \frac{\psi_{i_0-1,j_0} + \psi_{i_0,j_0-1}}{2}
    \\    &\le \frac{\sum_{l=2}^{i_0-1}   \psi_{l,1}  \xi\lrp{i_0-1,j_0,l} 
        + \sum_{l=2}^{j_0}  \psi_{1,l} \xi\lrp{i_0-1,j_0,l} }{2}
   \\ & + \frac{\sum_{l=2}^{i_0}  \psi_{l,1} \xi\lrp{i_0,j_0-1,l}  + \sum_{l=2}^{j_0-1}   \psi_{1,l} \xi\lrp{i_0,j_0-1,l} }{2} 
  \\  & = \sum_{l=2}^{i_0} \psi_{l,1}\xi\lrp{i_0,j_0,l} + \sum_{l=2}^{j_0}  \psi_{1,l}  \xi\lrp{i_0,j_0,l}.
\end{split}
\end{align*}
Thus, the induction hypothesis is satisfied.
    Therefore, for any $\bi=(i,j)$,
    \begin{equation*}
        \psi_{\bi} \le \sum_{l=1}^i \psi_{l,1} \xi\lrp{i,j,l} + \sum_{l=1}^j \psi_{1,l} \xi\lrp{i,j,l},
    \end{equation*}
    which leads to
\begin{equation*}
    \tau_{k,\bi}
    \le i+j + \sum_{l=1}^i \tau_{k,\lrp{l,1}} \xi\lrp{i,j,l}+ \sum_{l=1}^j \tau_{k,\lrp{1,l}} \xi\lrp{i,j,l}.
\end{equation*}
    Then, the following inequality holds:
    \begin{equation}\label{eq:2d_combinatorics}
        \tau_{k,\bi} - cN -b \le 
        \sum_{l=1}^i \lrp{\tau_{k,\lrp{l,1}} - cN - b} \xi\lrp{i,j,l}  +\sum_{l=1}^j \lrp{\tau_{k,\lrp{1,l}} - cN-b} \xi\lrp{i,j,l} + i + j.
    \end{equation}
    Define $\Delta_{k,\bi} = \tau_{k,\bi} - cN - b.$
For $\bi\in I_0$, without loss of generality, consider the case $j=N-1$.
   The second term of Eq (\ref{eq:2d_combinatorics}) can be decomposed as
   \begin{align*}
  \begin{split}
       \sum_{l=1}^j \Delta_{k,\lrp{1,l}} \xi\lrp{i,j,l}
      & \le  \lrp{\sum_{l=1}^m + \sum_{l=m+1}^{2^n-1}}  \Delta_{k,\lrp{1,l}} \xi\lrp{i,j,l}\\ 
     & \le 1 + \sum_{l=1}^m \Delta_{k,\lrp{1,l}} \xi\lrp{i,j,l}
     \le 1 + \sum_{l=1}^m \lrp{c_1N + c_2} \xi\lrp{i,j,l}.
   \end{split}
   \end{align*}
   Consider two cases: $i\ge m^3+2$ and $i<m^3+2$. 
   If $i\ge m^3 + 2$, by \cref{lemma:stirling},
   \begin{align*}
   \begin{split}
       2\xi\lrp{i,j,l} \le \frac{1}{2\pi \sqrt{i-2}}\sqrt{1+ \frac{i- 2}{N -1 -l}}
       \le \frac{1}{2\pi m\sqrt{m}}\sqrt{1+ \frac{m^3}{N - m^3 -m}}.
       \end{split}
   \end{align*}
   Therefore, 
   \begin{align*}
   \begin{split}
       1 + \sum_{l=1}^m \lrp{c_1 N+c_2} \xi\lrp{i,j,l} 
       \le 1 + \frac{c_1N+c_2}{4\pi \sqrt{m}}\sqrt{1+ \frac{m^3}{N - m^3 -m}} = o\lrp{N},
   \end{split}
   \end{align*}
   where the last equality holds as $m = \cO\lrp{\log N}=\cO\lrp{n}$.
   If $i < m^3 +2$, 
   \begin{align*}
  \begin{split}
       \xi\lrp{i,j,l} \le  2^{ -N +1}\lrp{i + j-l-2}^{i-2 } \le 2^{ -N +1} \lrp{2N}^{m^3} \le 2^{ - N +1 + m^3 \lrp{\log_2 {N}+1}}.
   \end{split}
   \end{align*}
   Thus, 
  \begin{align*}
   \begin{split}
        1 + \sum_{l=1}^m \lrp{c_1 N+c_2} \xi\lrp{i,j,l} \le 1 + m\lrp{c_1N+c_2} 2^{ - N +1 + m^3 \lrp{\log_2 {N}+1}} = o\lrp{N}.
   \end{split}
   \end{align*}
   In both cases,
   \begin{equation*}
         \sum_{l=1}^j \Delta_{k,\lrp{1,l}} \xi\lrp{i,j,l} = o\lrp{N}.
   \end{equation*}
   Similarly, the first term can be decomposed as follows:
    \begin{align*}
  \begin{split}
    \sum_{l=1}^i \Delta_{\lrp{l,1}} \xi\lrp{i,j,l}
   & = \lrp{\sum_{l=1}^m + \sum_{l=m+1}^{2^n-1}}     \Delta_{\lrp{l,1}} \xi\lrp{i,j,l} \\
   & \le 1 + \sum_{l=1}^m \Delta_{k,\lrp{l,1}} \xi\lrp{i,j,l}
     \le 1 +  \sum_{l=1}^m \lrp{c_1 N+c_2}\xi\lrp{i,j,l}.
  \end{split}
  \end{align*}
Again consider two cases: $i \ge m^3 +2 $ and $i < m^3 +2 $.
     If $i \ge m^3 +2 $, then, by \cref{lemma:stirling},
    \begin{align*}
  \begin{split}
     & \sum_{l=1}^i \lrp{\tau_{k,\lrp{l,1}} - cN - b} \xi\lrp{i,j,l} \\
     & \le 
         1 + \sum_{l=1}^m \Delta_{k,\lrp{l,1}} 2^{-i-j+l+1} \frac{\sqrt{i+j-l-2}}{2\pi \sqrt{\lrp{j-2}\lrp{i-l}}} 2^{i+j-l-2}\\
    & \le 1 + \sum_{l=1}^m \Delta_{k,\lrp{l,1}}\frac{1}{2\pi\sqrt{i-l}}  \\ 
     &\le 1 + \sum_{l=1}^m \lrp{c_1N + c_2}\frac{1}{2\pi\sqrt{m^3-m}}  = o\lrp{N}.
    \end{split}
    \end{align*}
    If $i<m^3+2$,
    \begin{align*}
  \begin{split}
         \xi\lrp{i,j,l} \le 2^{-N+1}\lrp{i+j-l-2}^{i-l} \le 2^{-N+1}\lrp{2N}^{m^3+2}
         \le 2^{-N+1 + \lrp{\log_2 N +1}\lrp{m^3+2}}.
    \end{split}
    \end{align*}
    Therefore, 
    \begin{align*}
  \begin{split}
        1 +  \sum_{l=1}^m \lrp{c_1 N+c_2}\xi\lrp{i,j,l} \le 1 +m\lrp{c_3 N+b_1} 2^{-N+1 + \lrp{\log_2 N +1}\lrp{m^3+2}} = o\lrp{N}.
    \end{split}
    \end{align*}
    We conclude that 
    \begin{equation*}
          \tau_{k,\bi}  \le cN +b + i+j + o\lrp{N} =c N +i +j +o\lrp{N},
    \end{equation*}
Now, for a general variable coefficient function $\coe$, to distinguish $\tau_{k,\bi}$ from the case of the constant function $\coe$ case, we denote it as $\tau^{\coe}_{k,\bi}$.
We observe that,
\begin{equation*}
    \left|\frac{\coe_{C_{k,\bi \pm q_j}}}{\coe_{C_{k,\bi}}} - 1\right| 
    \le \frac{h}{2^{n + \mathcal{L}\lrp{C_k}}}\frac{\max_{j\in \{1,2\}} \frac{\partial \coe}{\partial x_j}}{\coe_{C_{k,\bi}}}+\cO\lrp{2^{-2n}}.
\end{equation*}
Thus, 
\begin{equation*}
     \left|\frac{A_{C_{k,\bi},C_{k,\bi \pm q_j}}}{\coe_{C_{k,\bi}}} -1\right| = \left|\frac{\frac{2}{\frac{1}{\coe_{C_{k,\bi}}} + \frac{1}{\coe_{C_{k,\bi \pm q_j}}}}}{\coe_{C_{k,\bi}}} -1\right|
     \le {\frac{h}{2^{n +1+ \mathcal{L}\lrp{C_k}}}\frac{\lrp{\max_{j}\frac{\partial \coe}{\partial x_j}}_{C_{k,\bi}}}{\coe_{C_{k,\bi}}}}
     +\cO\lrp{2^{-2n}}.
\end{equation*}
By \cref{lemma:tau_evolove}, for $\bi = \lrp{i, j}$ such that $i,j\ge 2$, the following relation holds:
 \begin{align*}
  \begin{split}
    \tau^{\coe}_{k,\bi} &\le 1 + \frac{A_{C_{k,\bi},C_{k,\bi+q_1}} + A_{C_{k,\bi},C_{k,\bi+q_2}} } {\frac{A_{C_{k,\bi},C_{k,\bi-q_1}}}{\tau^{\coe}_{k, \bi-q_1}} + \frac{A_{C_{k,\bi},C_{k,\bi-q_2}}}{\tau^{\coe}_{k, \bi-q_2}}}\\
    & \le 1 + \lrp{1 + {\frac{h}{2^{n + \mathcal{L}\lrp{C_k}}}\frac{\lrp{\|\nabla \coe\|}_{C_{k,\bi}}}{\coe_{C_{k,\bi}}}}+\cO\lrp{2^{-2n}}}\frac{2}{\frac{1}{\tau^{\coe}_{k,\bi-q_1}} + \frac{1}{\tau^{\coe}_{k,\bi-q_2}}}.
\end{split}
\end{align*}
We can easily check that the following inequality between $\tau$ and $\tau^{\coe}$ defined using constant and variable coefficient function, respectively, holds:
 \begin{align*}
  \begin{split}
    \tau^{\coe}_{k,\bi} &\le \lrp{1 + {\frac{h}{2^{n + \mathcal{L}\lrp{C_k}}}\frac{\lrp{\|\nabla \coe\|}_{C_{k,\bi}}}{\coe_{C_{k,\bi}}}}+\cO\lrp{2^{-2n}}}^{i+j} \tau_{k,\bi}\\
    & \le \exp\lrp{ {\frac{h(i+j)}{2^{\mathcal{L}\lrp{C_k}}+n} \sup\frac{\|\nabla \coe\|}{\coe} } +\cO\lrp{2^{-n}}}\tau_{k,\bi}
    ,
\end{split}
\end{align*}
    and this completes the proof.
\end{proof}

\begin{lemma}\label{lemma:stirling}
    For $i, j\in \mathbb{N}$, the following inequality holds:
        \begin{equation*}
             2^{-i-j}\binom{i+j}{i} \le \frac{\sqrt{i+j}}{2\pi \sqrt{ij}}.
        \end{equation*}
    \end{lemma}
    \begin{proof}
            By Stirling's approximation, the following inequality holds:
 \begin{equation*}
     \binom{i+j}{i}\le \frac{{e^{-i-j+1}\lrp{{i+j}}^{i+j+0.5}}}{2\pi e^{-i} i^{i+0.5} \times e^{-j} j^{j+0.5}}
     \le \frac{\lrp{{i+j}}^{i+j+0.5}}{2\pi i^{i+0.5}j^{j+0.5}} 
     \le \frac{\sqrt{i+j}}{2\pi \sqrt{ij}} 2^{i+j}
     ,
 \end{equation*}
The last inequality follows from the Jensen's inequality:
 \begin{equation*}
     \frac{i}{i+j}\log\lrp{i}+\frac{j}{i+j}\log\lrp{j}
     \ge \log\lrp{ \frac{i^2+j^2}{i+j}}\ge \log\lrp{\frac{i+j}{2}}. 
 \end{equation*}
 Therefore, we obatin
 \begin{equation*}
     2^{-i-j}\binom{i+j}{i} \le  \frac{\sqrt{i+j}}{2\pi \sqrt{ij}} .
 \end{equation*}
 This completes the proof.
    \end{proof}

\subsection{Proof of \cref{theorem:3d_octree_main}}\label{subsec:proof:3d_octree_main}
\begin{proof}
    We first prove that there exists a constant $C$ such that for any $C_k\in O_0$, if the following inequality holds for any $\bi=\lrp{i_1,i_2,i_3}$ such that $i_1,i_2,i_3 < 2^n$,
    \begin{equation*}
        \tau_{k,\bi} \le c \cdot 2^{n},
    \end{equation*}
    then, for any $\bi\in I$ and $n\in \mathbb{N}$,
    \begin{equation*}
        \tau_{k,\bi}\le  cC \cdot 2^{n}.
    \end{equation*}
    If $i_1 = 2^n$ and $i_2,i_3<2^n$, we have
    \begin{equation*}
        \tau_{k,\bi}\le \frac{2\cdot 2^{-\mathcal{L}\lrp{C_k}} + 2^{-\mathcal{L}\lrp{C_k} +1}}{\frac{2^{-\mathcal{L}\lrp{C_k}}}{\tau_{k,\bi-q_2}} + \frac{2^{-\mathcal{L}\lrp{C_k}}}{\tau_{k,\bi-q_3}}}\lrp{1+\frac{P_{\coe}}{2^n}}
        \le 2c\lrp{1+\frac{P_{\coe}}{2^n}}\cdot 2^n.
    \end{equation*}
    Similar results hold for the cases $i_2 = 2^n$ and $i_1,i_3<2^n$, and $i_3 = 2^n$ and $i_1,i_2<2^n$.
    If $i_1 = i_2 = 2^n$ and $i_3<2^n$, we have
    \begin{equation*}
        \tau_{k,\bi}\le \frac{2^{-\mathcal{L}\lrp{C_k}} + 2\cdot 2^{-\mathcal{L}\lrp{C_k} +1}}{\frac{2^{-\mathcal{L}\lrp{C_k}}}{\tau_{k,\bi-q_1}} + \frac{2^{-\mathcal{L}\lrp{C_k}}}{\tau_{k,\bi-q_2}}} \lrp{1+\frac{P_{\coe}}{2^n}}
        \le 5c\lrp{1+\frac{P_{\coe}}{2^n}}^2\cdot 2^n.
    \end{equation*}
     Similar results holds for the cases $i_2 = i_3 = 2^n$ and $i_1 = i_3 = 2^n$.
     If $i_1 = i_2 = i_3= 2^n$,
         \begin{equation*}
        \tau_{k,\bi}\le \frac{3 \cdot 2^{-\mathcal{L}\lrp{C_k} +1}}{\frac{2^{-\mathcal{L}\lrp{C_k}}}{\tau_{k,\bi-q_1}} + \frac{2^{-\mathcal{L}\lrp{C_k}}}{\tau_{k,\bi-q_2}} + \frac{2^{-\mathcal{L}\lrp{C_k}}}{\tau_{k,\bi-q_3}}}\lrp{1+\frac{P_{\coe}}{2^n}}
        \le 10c\lrp{1+\frac{P_{\coe}}{2^n}}^3\cdot 2^n.
    \end{equation*}
    Therefore, the constant $C$ exists.

    Now, we use mathematical induction to prove the following statement:
    For any $C_k\in O_0$ there exists a constant $c_k$ such that for any $\bi=\lrp{i_1,i_2,i_3}$ such that $i_1,i_2,i_3 <2^n$ and $n\in \mathbb{N}$, the following inequality holds:
    \begin{equation*}
        \tau_{k,\bi}\le c_k 2^n.
    \end{equation*}
    First, consider $C_k\in O_0$ such that 
    \begin{equation*}
        C_k \nprec_0 C_{k'}, 
    \end{equation*}
    for any $C_{k'}$.
    It is obvious that for $\bi =\lrp{i_1,i_2,i_3}$ such that $i_1,i_2,i_3\le 2^n$, 
    \begin{equation*}
        \tau_{k,\bi}\le  \exp\lrp{ {\frac{3h}{2^{\mathcal{L}\lrp{C_k}}} \sup\frac{\|\nabla \coe\|_1}{\coe} }} 3\cdot 2^{n}.
    \end{equation*}
  Now assume that the induction hypothesis holds for any $C_k$ such that $C_{k}\prec_O C_{k_0}$.
That is, there exists a constant $c_k$ such that the following inequality holds for any $n\in \mathbb{N}$ and $\bi=\lrp{i_1,i_2,i_3}$ such that $i_1, i_2, i_3 <2^n$:
  \begin{equation*}
      \tau_{k,\bi}\le c_k 2^n.
  \end{equation*}
  Then, for any $\bi\in I$,
    \begin{equation*}
      \tau_{k,\bi}\le Cc_k 2^n.
  \end{equation*}
  Then, for $\bi = \lrp{1,i_2,i_3}$, there exists a cell $C_{\lrp{k, \lrp{2^n, i_2', i_3'}}}$ such that $C_{\lrp{k, \lrp{2^n, i_2', i_3'}}}\prec_O C_{\lrp{k_0, \lrp{1, i_2, i_3}}}$ and 
  \begin{equation*}
      \tau_{k_0,\bi}\le \frac{3 \cdot 2^{-\mathcal{L}\lrp{C_{k_0}}}}{\frac{2^{-\mathcal{L}\lrp{C_k}}}{\tau_{k, \lrp{2^n, i_2', i_3'}}} }\lrp{1+\frac{P_{\coe}}{2^n}}
      \le 6Cc_k\lrp{1+\frac{P_{\coe}}{2^n}} 2^n.
  \end{equation*}
  Similar results hold for $\bi = \lrp{i_1,1,i_3}$ and $\bi = \lrp{i_1,i_2,1}$.
  Then, for any $\bi= \lrp{i_1,i_2,i_3}$ such that $i_1,i_2,i_3 < 2^n$,
  \begin{equation*}
      \tau_{k_0,\bi}
      \le \max_{\substack{\bi' = \lrp{i_1,i_2,i_3},\\ i_1=1, i_2=1,\text{ or }i_3=1}} \tau_{k_0,\bi'} \lrp{1+\frac{P_{\coe}}{2^n}}^{2^n} \le \lrp{6Cc_k +3}e^{P_{\sigma}} 2^n.
  \end{equation*}
  Thus, the induction hypothesis is satisfied.
  This completes the proof.
\end{proof}


\bibliographystyle{siamplain}
\bibliography{mybib}
\end{document}